\documentclass[12pt]{article}
\usepackage{amssymb}
\usepackage{amsthm}
\usepackage{amsmath} 
\usepackage{xcolor}
\usepackage{float}
\usepackage{pifont}
\usepackage{tikz}
\usepackage{pgfplots}
\pgfplotsset{compat=1.16}
\usetikzlibrary{patterns}
\usetikzlibrary{arrows}
\usepackage{stmaryrd,bm}
\usepackage{soul}
\usepackage{lineno}
\usepackage{mathrsfs}
\usepackage{textgreek}
\usepackage{enumitem}
\usepackage{isomath}
\usepackage{rotating}
\usepackage{cite}
\usepackage{accents}
\usepackage[unicode,linkcolor=blue,colorlinks=true,citecolor=red]{hyperref}

\newcommand{\bb}{\bm b}
\newcommand{\dst}{\displaystyle}
\renewcommand{\leq}{\leqslant}
\renewcommand{\geq}{\geqslant}
\newcommand{\noin}{\noindent}
\newcommand{\la}{\langle}
\newcommand{\ra}{\rangle}

\newcommand*{\dt}[1]{%
  \accentset{\mbox{\large\bfseries .}}{#1}}
\newcommand{\bchi}{\bm \chi}
\newcommand{\ep}{\sigma}
\renewcommand{\rho}{\varrho}
\newcommand{\rhop}{\rho_{+}}

\renewcommand{\phi}{\varphi}
\newcommand{\kp}{k_{+}}
\newcommand{\km}{k_{-}}
\newcommand{\m}{\bm p}
\newcommand{\mb}{\bm m}

\newcommand{\br}{\bm x}

\newcommand{\qh}{\what{q}}
\newcommand{\bk}{{\bm k}}

\newcommand{\n}{\bm n}
\newcommand{\de}{\partial}
\newcommand{\vr}{\varrho}
\newcommand{\x}{{\bm x}}
\newcommand{\si}{{\bm \xi}}
\newcommand{\et}{{\bm \eta}}

\newcommand{\y}{{\bm y}}

\newcommand{\hpsi}{\psi}

\newcommand{\rf}[1]{(\ref{#1})}

\newcommand{\ut}{\tilde{u}}

\newcommand{\uh}{\what{u}}

\newcommand{\vt}{\tilde{v}}
\newcommand{\wh}{\what{w}}
\newcommand{\wt}{\tilde{w}}
\newcommand{\walpha}{\what{\alpha}}

\renewcommand{\setminus}{\smallsetminus}
\renewcommand{\tilde}{\widetilde}
\newcommand{\what}{\widehat}
\newcommand{\nonu}{\nonumber}
\newcommand{\oh}{\frac{1}{2}}

\newcommand{\dotx}{\dt{\x}}
\newcommand{\dotk}{\dt{\bk}}

\newcommand{\tr}{\mbox{\sf tr}}
\newcommand{\z}{{\bm 0}}

\newcommand{\A}{\mathsfbfit A}
\newcommand{\Am}{\bm A}
\newcommand{\B}{\mathsfbfit B}

\newcommand{\Csbi}{\mathsfbfit C}
\newcommand{\Dsbi}{\mathsfbfit D}

\renewcommand{\U}{{\mathsfbfit S}}
\renewcommand{\G}{\mathsfbfit G}

\newcommand{\Cm}{\bm C}
\newcommand{\Dm}{\bm D}
\newcommand{\Ho}{H^\oh}
\newcommand{\Ht}{H^\frac{3}{2}}

\newcommand{\M}{\bm M}
\newcommand{\Msbi}{\mathsfbfit M}
\newcommand{\Om}{\Omega}
\newcommand{\Omh}{\what{\Omega}}

\renewcommand{\P}{\mathsfbfit P}
\newcommand{\R}{{\mathbb R}}

\newcommand{\T}{{\mathsfbfit T}}
\newcommand{\F}{{\mathsfbfit F}}
\newcommand{\Id}{{\sf I}}
\renewcommand{\Im}{\bm I}
\newcommand{\J}{J}
\newcommand{\Q}{\bm Q}
\renewcommand{\O}{{\cal O}}

\newcommand{\N}{\mathsfbfit N}
\newcommand{\Np}{{\N}^{+}}
\newcommand{\Nm}{{\N}^{-}}
\newcommand{\Z}{{\mathbb Z}}

\newcommand{\vre}{\vr_{+}}
\newcommand{\gi}{\gamma_{-}}
\renewcommand{\ge}{\gamma_{+}}
\newcommand{\om}{\omega}
\newcommand{\s}{\Pi}
\newcommand{\ei}{\varepsilon}

\DeclareMathAlphabet{\mathpzc}{OT1}{pzc}{m}{it}

\date{}

\newtheorem{theorem}{Theorem}

\newtheorem{lemma}{Lemma}
\newtheorem*{definition}{Definition}

\newtheorem*{remark}{Remark}

\renewcommand{\theequation}{\thesection.\arabic{equation}}
\addtocounter{figure}{0}

\newcommand{\I}{\ensuremath{\mathrm{i}}}
\newcommand{\E}{\ensuremath{\mathrm{e}}}
\newcommand{\D}{\ensuremath{\mathrm{d}}}

\begin{document}

\title{Clusters of Bloch waves in three-dimensional periodic media}
\author{
Yuri~A. Godin\thanks{Email: ygodin@uncc.edu} and Boris Vainberg\thanks{Email: brvainbe@uncc.edu} \\
The University of North Carolina at Charlotte, \\
Charlotte, NC 28223 USA
}

\maketitle

\subsection*{Abstract}
We consider acoustic wave propagation through a periodic array of the inclusions of arbitrary shape. The inclusion size is much smaller than the array period while the wavelength is fixed.
We derive and rigorously justify the dispersion relation for general frequencies and show that there are exceptional frequencies for which the solution is a cluster of waves propagating in different directions with different frequencies so that the dispersion relation cannot be defined uniquely. Examples are provided for the spherical inclusions.

\bigskip
\noin
{\bf Keywords:}\\
Periodic media; phononic crystal; Bloch waves; dispersion relation; asymptotic expansion;
Dirichlet-to-Neumann operator;

\section{Introduction}
\setcounter{equation}{0}

The propagation of waves in periodic media reveals numerous effects of practical importance.
These include
bands and gaps in the wave spectrum, positive or negative group velocity, slowing down considerably the speed of light, nonreciprocal media, the self-collimation effect, wave filters, mirrors, and more that lead to the development of new devices (see \cite{JJWM:11} and references therein).

Dispersion of waves in such periodic structures can be studied numerically using a plane wave expansion, finite-difference or finite-element methods \cite{Kuchment:99}, or boundary-element methods \cite{JJWM:11}.
Deriving an explicit dispersion relation for the Floquet-Bloch waves in two and three-dimensional periodic media is an arduous problem and is usually performed numerically \cite{JJWM:11}. However, assuming that the wavelength is long compared to the period of the lattice or a characteristic size of the scatterers one can obtain an asymptotic approximation \cite{Kuchment:99a}.
One of the popular techniques to study dispersion relations in periodic media with small circular inclusions is matched asymptotic expansions which were applied to Dirichlet \cite{McIver:2009}, \cite{Craster:2017}
or Neumann \cite{McIver:07}, \cite{Craster:20} scatterers. The method was applied to the elastic waves propagating through a lattice of cylindrical cavities \cite{McIver:2011} as well as acoustic waves in three-dimensional lattices of arbitrary shape scatterers with the Neumann boundary condition \cite{Guo:14}.
A semi-analytical approach using the multipole expansion method is described in \cite{MMP:02}.
Explicit formulas for the effective dielectric tensor and the dispersion relation are obtained in \cite{GV:19} assuming that the cell size is small compared to the wavelength, but large compared to the size of the inclusions. Some other approaches are presented in \cite{Zalipaev:02,Craster:10,Vanel:17,Cherednichenko:06,Smyshlyaev:08,Joyce:17,Guzina:21}.

The present paper deals with the propagation of acoustic waves in an infinite medium containing a periodic array of identical inclusions of arbitrary shape with transmission conditions on their interfaces. The results can be immediately applied in the case of the Neumann boundary condition and the approach can be extended to the Dirichlet boundary condition.

Most of the previous studies constructed the functions that approximately satisfy the equation and the boundary conditions. Without an estimate on the inverse operator, one cannot guarantee that the constructed functions are close to the exact solution. We provide a rigorous justification that our approximation is asymptotically close to the exact solution. The rigorous approach reveals the existence of exceptional wave vectors for which the solution of the problem has the form of a cluster of waves propagating in different directions or with different frequencies. For those exceptional wave vectors, there is no single wave propagating in one direction with a particular spatial frequency, so the dispersion relation cannot be defined uniquely. This effect has some similarities with the Bragg reflection \cite{Brillouin:56} and refraction of waves.

The presence of small inclusions makes the inclusionless problem singularly perturbed. We avoid this formidable difficulty by reducing the problem to an operator equation on the surface of the sphere of fixed radius $R$ enclosing the inclusion. That allows us to use the standard perturbation theory and find the solution of the auxiliary problem in the form of a power series and construct rigorously the solution of the original problem.

The paper is organized as follows.
In Section \ref{form} we formulate the problem, introduce the notion of the exceptional Bloch vector and briefly describe the main results of the paper.
In Section \ref{out} we construct the inner and outer Dirichlet-to-Neumann (DtN) operators and reduce the singularly perturbed problem in question to a regular one for the perturbation of the zero eigenvalue of the difference of the DtN operators. Next, in Sections \ref{exter} and \ref{inner} we expand the DtN operators in the power series in terms of the small parameters of the problem.
The paper's main theorems on the structure of the solution of the problem and the dispersion for regular and exceptional Bloch vectors are formulated in Section \ref{main}.
These theorems are proven in Section \ref{proof}.
We illustrate our results by an example of a simple cubic lattice of spherical inclusions in Section \ref{spher}. Conclusions (Section \ref{concl}) contains also a discussion of spectral gaps.
The proof of the expansion of the inner DtN operator is somewhat tedious and is relegated to Appendix.

\section{Formulation of the problem and description of results}
\label{form}
\setcounter{equation}{0}

We consider the propagation of acoustic waves through an infinite medium containing a periodic array
of small identical inclusions.
\begin{figure}[th]
\begin{center}
\begin{tikzpicture}[scale=1.1,>=triangle 45]
 \begin{scope}[x={(4cm,0cm)},y={({cos(30)*2.5cm},{sin(30)*2.5cm})},
    z={({cos(70)*3cm},{sin(70)*3cm})},line join=round,fill opacity=0.5]
  \draw[fill=lightgray!50] (0,0,0) -- (0,0,1) -- (0,1,1) -- (0,1,0) -- cycle;
  \draw[fill=lightgray!50] (0,0,0) -- (1,0,0) -- (1,1,0) -- (0,1,0) -- cycle;
  \draw[fill=lightgray!50] (0,1,0) -- (1,1,0) -- (1,1,1) -- (0,1,1) -- cycle;
  \draw[fill=lightgray!50] (1,0,0) -- (1,0,1) -- (1,1,1) -- (1,1,0) -- cycle;
  \draw[fill=lightgray!50] (0,0,1) -- (1,0,1) -- (1,1,1) -- (0,1,1) -- cycle;
  \draw[fill=lightgray!50] (0,0,0) -- (1,0,0) -- (1,0,1) -- (0,0,1) -- cycle;
  \draw [->,very thick] (0,1,0) -- (1,1,0);
  \draw [->,very thick] (0,1,0) -- (0,0,0);
  \draw [->,very thick] (0,1,0) -- (0,1,1);
 \end{scope}
 \draw [->, thick] (4.15,1.55) -- (5.05,1.55) node [above right] {\Large  $x_2$};
 \draw [->, thick] (2.85,1.05) -- (2.35,0.35) node [left] {\Large  $x_1$};
 \draw [->, thick] (3.35,2.55) -- (3.35,3.35) node [above right] {\Large  $x_3$};
 \shade[ball color = lightgray, opacity = 0.9] (4.5,2.9,3) circle (1cm);
 \def\eggheight{3mm}
 \begin{turn}{45}
  \path[ball color=orange!80!gray]
  plot[domain=-pi:pi,samples=100,shift={(3.58,-1.15)}]
  ({.78*\eggheight *cos(\x/4 r)*sin(\x r)},{-\eggheight*(cos(\x r))})
  -- cycle;
  \end{turn}
  \node [above] at (3.35,1.45) {$\Om$};
  \node [above] at (5.8,4.0,2) {\Large$\Pi$};
  \node [above] at (3.1,1.95) {\large$B_R$};
  \node [above] at (6.4,0.5) {\Large $\bm \ell_2$};
  \node [above] at (-0.3,0) {\Large $\bm \ell_1$};
  \node [above] at (3,4.1) {\Large $\bm \ell_3$};
\end{tikzpicture}
\end{center}
\caption{The cell of periodicity $\Pi$ containing a ball $B_R$ of radius $R$ which encloses an inclusion $\Om$.}
\label{cell}
\end{figure}
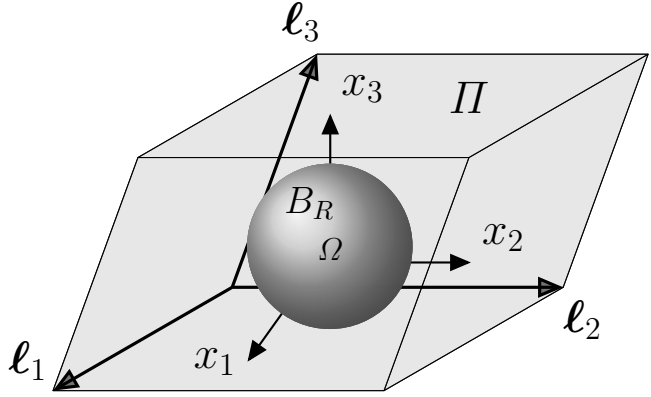
The periodicity of the medium is defined by the vectors ${\bm \ell}_1, {\bm\ell}_2, {\bm\ell}_3$. We fix the cell of periodicity $\Pi$ (a parallelepiped) in such a way that $\de\Pi$ does not intersect the inclusion. Denote by $\Om$ the domain occupied by the inclusion in $\Pi$ (see Figure \ref{cell}). We choose the origin of the coordinate system in $\Om$ and assume that $\Om$ is small with the size of order $a>0$. More precisely,  $\Om = \Om (a)$ is obtained from an $a$-independent domain $\what{\Om}$ by the contraction with the coefficient $a^{-1}$, i.e. the transformation $\x \to a\si $ maps $\Om (a) \subset \R^3_{\x}$ into $\Omh \subset \R^3_{\si}$.
We assume that $ \de\Om  \in C^{1,\beta}$, i.e. the functions describing the boundary have first derivatives that belong to the H\"{o}lder space with index $\beta$.

Propagation of acoustic time-harmonic waves with time frequency $\om$ is governed by the equation
\begin{equation}
 \Delta u + k_{\pm}^2  u = 0,
 \label{eq}
\end{equation}
where $u$ is the amplitude of the excess pressure and
$k_{\pm}$ are the wave numbers in the medium and the inclusion, respectively. Here and on the subscript $\pm$ refers to the value of the quantity outside/inside of the inclusion. We assume that the mass density is piecewise-constant: $\rho (\x) = \rho_{\pm}$.
The wave numbers $k_\pm$ are related to the frequency $\om$ by $k_{\pm} = \om/c_{\pm}$, where
$c_\pm = 1/\sqrt{\gamma_\pm\rho_\pm}$ is the speed of the wave propagation and $\gamma_\pm$ is the adiabatic bulk compressibility modulus.
We assume that the inclusions are penetrable and therefore impose the transmission conditions
on their boundaries
\begin{align}
\label{bc1}
\left. \left\llbracket u (\br) \right\rrbracket \right . &=0, \\[2mm]
\left. \left\llbracket \frac{1}{\vr(\br)} \frac{\de u (\br)}{\de \n} \right\rrbracket \right. &=0,
\label{bc2}
\end{align}
where $\vr(\br)$ is the mass density and $\n$ is the outward normal vector. The brackets $\llbracket \cdot \rrbracket$ denote the jump of the enclosed quantity across the boundaries of the inclusions.

We study the propagation of Floquet-Bloch waves that satisfy \rf{eq}-\rf{bc2} and have the form
\begin{equation}
 u(\br) = \Phi (\br) \,\E^{-\I \bk \cdot \br},
\end{equation}
where $\bk = (k_1,k_2,k_3)$ is the wave vector and $\Phi (\br)$ is a periodic function
with the periods of the lattice. The latter condition is equivalent to
\begin{equation}
 \rrbracket  \E^{\I\bk \cdot \br} u(\br) \llbracket =0,
\label{FB1}
\end{equation}
where the inverted brackets $\rrbracket  \cdot \llbracket$ denote the jump of the enclosed expression and their first derivatives across the opposite sides of the cells of periodicity.

Thus, we arrive at the following problem for the amplitude $u$ in the fundamental cell of periodicity $\s$:
\begin{equation}
\begin{array}{l}
\Delta u_{-} + k^2_{-}  u_{-} = 0, \quad u_{-} \in H^2(\Om), \\[2mm]
\Delta u_{+} + k^2_{+}  u_{+} = 0, \quad u_{+} \in H^2 (\s\setminus\Om),
\end{array}
\label{Hz1}
\end{equation}
\vspace{-3ex}
\begin{align}
\label{bc1a}
\left\llbracket u (\br) \right\rrbracket =0, \quad
\left\llbracket \frac{1}{\vr(\br)} \frac{\de u (\br)}{\de \n} \right\rrbracket =0, \quad
\rrbracket  \E^{\I\bk \cdot \br} u(\br) \llbracket =0.
\end{align}
Here $H^2(\Om)$ and $H^2 (\s \setminus \Om)$ are the Sobolev spaces.
We look for the dispersion relation, i.e. the relation between $\om=\kp c_{+}$ and the wave vector $\bk$  for which the problem  \rf{Hz1}-\rf{bc1a} has a nontrivial solution.

In the inclusionless case, there is a simple dispersion relation between the time frequency $\om$ and the spatial frequency $|\bk|$. Namely, the unperturbed problem with $\kp = |\bk|$ has the solution $\E^{-\I \bk \cdot \x}$, and therefore the dispersion relation is $\om=c|\bk|$, where $c=c_+ = 1/\sqrt{\ge\vre}$ is the speed of waves in the host medium. In the presence of inclusions the dispersion relation $\om=H(\bk,a)$  is more complicated and our goal is to find it when $a$ is small.

In the standard fashion, we introduce the basis vectors $\bb_1, \bb_2, \bb_3$ of the reciprocal lattice $\mathbb Z^3_b$
\begin{align}
 \label{}
 {\bm \ell}_i \cdot \bb_j = 2\pi \delta_{i,j},
\end{align}
where $\delta_{i,j}$ is the Kronecker delta. Points $\mb$ of the lattice are
\begin{align}
 \mb = m_1 \bb_1 + m_2 \bb_2 + m_3 \bb_3,
 \label{mb}
\end{align}
where $m_1, m_2, m_3$ are integers. If $\Pi$ is a cube $[-\pi,\pi]^3$, then $\mathbb Z^3_b$ is the standard lattice $\mathbb Z^3$ of the points $\mb = (m_1 , m_2,  m_3 )$.

It is important to keep in mind that the space of Bloch waves is not always one-dimensional even for the unperturbed problem.

\begin{definition}
 A point $\bk$ is called exceptional if there are non-trivial vectors $\mb = \mb_s, ~ 2 \leq s \leq n,$ of the form \rf{mb} such that $|\bk| = |\bk - \mb|$. The number $n \geq 2$ is called the order of the exceptional point. We set $\mb_1 = \z$.
\end{definition}
\begin{remark}
 Geometrically, $n$ is the number of  points of the reciprocal lattice $\mathbb Z^3_b$ (including the origin) on the sphere of radius $|\bk|$ centered at $\bk$.
\end{remark}

Consider the unperturbed problem
\begin{equation}
\label{unp}
\Delta u + |\bk|^2  u = 0, \quad u \in H^2 (\s), \quad
\rrbracket  \E^{\I\bk \cdot \br} u(\br) \llbracket =0.
\end{equation}
 \begin{lemma}\label{lex}
 The solution space of the problem \rf{unp} is $n$-dimensional, and it is spanned by functions $\psi_s=\E^{-\I (\bk -\mb_s)\cdot \x}, ~ 1\leq s\leq n,$ where $\mb_s $ are the points of $\mathbb Z^3_b$ such that $|\bk| = |\bk - \mb_s|$.
 \end{lemma}
 \begin{proof}
If $u$ is a solution of \rf{unp}, then function $v(\x)=\E^{\I \bk \cdot \x}u(\x)$ is a periodic solution of the equation $\Delta v -2i\bk\cdot \nabla v = 0,$
and $v$ can be extended in the Fourier series
\[
v(\x)=\sum_{\mb\in\mathbb Z_b^3}a_{\mb}\E^{-\I \mb\cdot \x}.
\]
We put the series into the equation and obtain $(2\bk \cdot \mb-|\mb|^2)a_{\mb}=0$. Thus $a_{\mb}$ can be different from zero if and only if $2\bk \cdot \mb-|\mb|^2=0$. The latter relation is equivalent to $|\bk| = |\bk - \mb|$. Hence $v(\x)=a_{\z}$ is a constant if $\bk$ is not exceptional, and $v(\x)$ is a linear combination of functions $\E^{\I \mb_s\cdot \x},~1\leq s\leq n,$ if $\bk$ is an exceptional vector.

 \end{proof}

While the space of Bloch waves is multidimensional when $\bk$ is an exceptional point, this fact is not important for the unperturbed problem since these waves are untangled and propagate in different directions. We will see that the situation is much more complicated when the problem has inclusions.

From the definition, it follows that $\bk$ is exceptional if it belongs to at least one of the planes $2\bk \cdot \mb=|\mb|^2,\mb\in\mathbb Z^3_b ,$ and it has order $n+1$ if it belongs to $n$ such planes. The distance from the origin to these planes goes to infinity as $|\mb|\to\infty$. One can easily see that the unperturbed problem with $\kp = |\bk|$ has a multidimensional space of Bloch waves if and only if $\bk$ is exceptional, and the dimension of this space is $n \geq 2$.

When a small perturbation is introduced, the Bloch wave with a non-exceptional $\bk$ changes slightly outside of the perturbed area. We will show that $\bk$ is a smooth function of $a$, i.e., the dispersion relation is well defined. It will be shown that the dispersion relation has the following  asymptotics:
\begin{align}
 \om=c_{+}|\bk| (1+\mu a^3+\O(a^4)), \quad a\to 0, \quad \mu = \mu(\dotk), \quad \dotk = \bk/|\bk|,
 \label{om_reg}
\end{align}
where $c_{+}$ is the speed of waves in the host medium and the coefficient $\mu$ is defined by the properties of the media and $\Omega$, see Theorem \ref{t1}. This Bloch wave propagates in the direction of vector $\bk$.

When an exceptional wave vector $\bk$ of order $n$ is fixed, the perturbation excites simultaneously $n$ Bloch waves
$u_s (\x)$ and each of them has its time frequency
\begin{align}
 \om_s=c_{+}|\bk| (1+\mu_s a^3+\O(a^4)), \quad a\to 0,  \quad \mu_s = \mu_s(\dotk).
 \label{2.12}
\end{align}
Moreover, each solution $u_s$
is a cluster of waves propagating in directions $\bk-\mb_j, 1\leq j\leq n$, see Theorem \ref{t3}.
One can obtain a wave propagating in one direction by taking a linear combination of solutions $u_s$, but the combination would contain the terms with different time frequencies and these terms are solutions of different equations \rf{Hz1} with $\kp = \om_s/c_{+}$.

Let us invert the relation between $\bk$ and $\om$ and find $\bk$ as a function of $\om$. Since $\om$ is a number and $\bk$ is a vector, we look for the relation
\begin{align}
|\bk| = f(\om, \dotk, a).
\label{2.13}
\end{align}
For the unperterbed problem, $\om = c_{+} |\bk|$. Thus we fix $\om$ and a vector $\bk^\ast$ with an arbitrary direction such that $\om = c_{+} |\bk^\ast|$, and consider Bloch vectors $\bk = c\,\bk^\ast$ with $|c-1| \ll 1$. If $\bk^\ast$ is not exceptional, then the relation \rf{2.13} can be found by solving \rf{om_reg} for $|\bk|$. This cannot be done for \rf{2.12}. Theorem \ref{t4} concerns the case when the time frequency $\omega$ is fixed in such a way that $k_+=\omega/c_+=|\bk^\ast|$, where $\bk^\ast$ is an exceptional point, and $\bk$  is close to $\bk^\ast$ and have the same direction. It is proven there that
the wave vector $\bk$ has $n$ values $\bk_s$ close to $\bk^\ast $:
\[
\bk_s = \bk^\ast (1 + \nu_s a^3 + \O(a^4))
\]
for which problem \rf{Hz1}-\rf{bc1a} with $k_+=|\bk^\ast|$ and $\bk=\bk_s$  has a non-trivial solution $u_s$. Each of them is a cluster of waves propagating in the directions $\bk_s-\mb_j, 1\leq j\leq n$. One can construct a linear combination of $u_s$ to obtain a wave propagating in one direction, for example, in the direction of $\bk^\ast$, but the combination would contain the terms with different spatial frequencies $|\bk_s|$. Thus, a dispersion relation is not defined uniquely for the exceptional wave vectors.

The results obtained in the paper are asymptotic. The appearance of clusters instead of a single wave is justified if the size of the inclusion is sufficiently small compared to the periods of the cell.

An example with circular inclusions is provided in Section \ref{spher}. The Appendix contains proof of Theorem \ref{t2} on the
asymptotics of the solution of an auxiliary problem in the ball containing the inclusion with the Dirichlet boundary condition on the sphere. This theorem has an idependent value since it provides the asymptotic expansion in $a$ of the solution of the Dirichlet problem outside of the inclusion and for the Dirichlet-to-Neumann map  without the construction of a tedious asymptotics near the inclusion.

\section{Outline of the approach}
\label{out}
\setcounter{equation}{0}
In this paper, we simplify our approach introduced in \cite{GV:20} for spherical inclusions and extend it to the inclusions of arbitrary shapes.
Note also that our results in \cite{GV:20} are valid only for non-exceptional $\bk$.

To find the dispersion relation we introduce function $\varepsilon=\varepsilon(a,\omega,\bk)$ such that $ \kp=\om/c_+ = (1+\varepsilon) |\bk|$ and then determine $\varepsilon$ for which the problem  \rf{Hz1}-\rf{bc1a} has a nontrivial solution. 
Since $a=0, \varepsilon=0$ for the inclusionless problem, we look for such function $\varepsilon$ that vanishes as $a \to 0$.
The exact asymptotic behavior of $\varepsilon$ will be found.

The problem \rf{Hz1}-\rf{bc1a} is a singular perturbation of the inclusionless problem, and its solution is rather complicated in a neighborhood of the inclusion. We derive the dispersion relation from an auxiliary regularly perturbed problem.
To this end, we enclose the inclusion $\Om = \Om(a)$ in a ball $B_R \subset \Pi$ of radius $R > a$ centered at the origin, split $\Pi$ into the ball $B_R = \{ |{\br} | < R\}$ and its complement $\Pi \setminus B_R$ (see Figure \ref{cell}), and consider two separate problems in  $ B_R$ and $\Pi \setminus B_R$:
 \begin{align}
 \left( \Delta + k^2_{\pm} \right) u_{\pm}(\br) &= 0, \quad \br\in B_R, \; \x \notin \de \Om, ~~ \left\llbracket u (\br) \right\rrbracket =  \dst
\left\llbracket \frac{1}{\vr(\br)} \frac{\de u (\br)}{\de \n} \right\rrbracket=0,~~ \left. u_{+}\right|_{r=R}=\psi.
 \label{uin}\\[2mm]
 \label{k-egv}
 \left( \Delta + k^2_{+} \right) v(\br) &= 0, \quad  \br\in \Pi \setminus B_R,~~  \left\rrbracket  \E^{\I \bk \cdot \br} v(\br) \right\llbracket =0, ~~ \left. v\right|_{r=R}=\psi, \quad \kp = (1+\varepsilon) |\bk|.
\end{align}
Here  $u_{-} \in H^2 (\Om)$, $u_{+} \in H^2 (B_R \setminus \Om)$, $v \in H^2 (\Pi \setminus B_R)$.

In the unperturbed case, $\varepsilon =0, ~ a=0$, problems \rf{uin}, \rf{k-egv} are uniquely solvable for all values of $R$, except possibly a discrete set $\{R_i\}$. We fix an $R\notin \{R_i\}$. It will be shown that solutions of \rf{uin}, \rf{k-egv} are still unique for small values of $a$ and $\varepsilon$. We define  operators $\Nm_{a,\varepsilon}$ and $\Np_{\bk,\varepsilon}$ as the Dirichlet-to-Neumann (DtN) operators for problems \rf{uin}, \rf{k-egv}  with the derivatives in the direction of $r$:
\begin{equation}
\label{}
\Nm_{a,\varepsilon}, \Np_{\bk,\varepsilon}:\Ht(\de B_R)\to \Ho(\de B_R),
\quad~~ \Nm_{a,\varepsilon}  \psi = \left.\frac{\de u}{\de r}\right|_{r=R}, \quad  \quad
\Np_{\bk,\varepsilon}  \psi = \left.\frac{\de v}{\de r}\right|_{r=R},
\end{equation}
where $H^{s}(\de B_R)$ is the  Sobolev space of functions on $\de B_R$.

We provide below two lemmas useful for the description of our approach and the explanation of its advantage.

\begin{lemma}\label{l2}
\leavevmode
\begin{enumerate}
 \item[(1)]
 Operators $\Np_{\bk,\varepsilon},~\Nm_{a,\varepsilon}$ and
 \begin{equation}\label{npm}
\Np_{\bk,\varepsilon}-\Nm_{a,\varepsilon}:\Ht(\de B_R)\to \Ho(\de B_R)
 \end{equation}
are Fredholm.
\item[(2)]
Operators $\Np_{\bk,\varepsilon}$ and $\Nm_{a,\varepsilon}$ are symmetric in $L^{2}(\de B_R)$. Thus, for example,
 \begin{equation}\label{sym}
  \int_{r=R} \left( \Np_{\bk,\varepsilon} \psi \right) \overline{\phi}\, \D S = \int_{r=R} \psi \, \overline{\left(\Np_{\bk,\varepsilon} \phi\right)}\, \D S,  \quad \psi,\phi\in \Ht(\de B_R).
 \end{equation}
\end{enumerate}
\end{lemma}
The proof is given in \cite{GV:20} (the shape of $\Om$ is irrelevant in the proof). The Fredholm property of the operators follows from their ellipticity (the ellipticity of the DtN map and its symbol can be found in \cite{vg}). The symmetry of the operators is a direct consequence of Green's formula.

\begin{lemma}\label{l1a}
\leavevmode
Relation $\psi=u|_{r=R}$ is a one-to-one correspondence between solutions $u$ of \rf{Hz1}-\rf{bc1a} and solutions $\psi\in \Ht(\de B_R)$ of
\begin{align}
\label{nn}
(\Np_{\bk,\varepsilon}-\Nm_{a,\varepsilon})\,\psi=0.
\end{align}
\end{lemma}
\begin{proof}
Let $u$ be a solution of  \rf{Hz1},\rf{bc1a}. Then the Dirichlet and Neumann data for $u_+$ coincide with the data for $u_-$, and therefore \rf{nn} holds. The inverse statement is a consequence of the fact that the extension of any solution of \rf{uin} onto domain $\Pi \setminus B_R$ using the solution of  \rf{k-egv} gives a solution of  \rf{Hz1},\rf{bc1a} provided that the Dirichlet and Neumann data on $\de B_R$ of the solution of \rf{uin} coincide with the corresponding data of the solution of \rf{k-egv}.
\end{proof}

Lemma \ref{l1a} is the key point in our paper. It reduces the dispersion relation to the relation between $\bk$ and $\om$  for which $\Np_{\bk,\varepsilon} - \Nm_{a,\varepsilon}$ has a zero eigenvalue. The advantage of our approach is based on the fact that $\Np_{\bk,\varepsilon}$ does not depend on the inclusion and $\Nm_{a,\varepsilon}$ is an infinitely smooth operator-function of $a$.
This allows us to find the dispersion relation using the standard perturbation theory. Technically, it is more convenient to replace $\Np_{\bk,\varepsilon}$ and $\Nm_{a,\varepsilon}$
by $\Np_{\bk,\varepsilon}- \Nm_{0,\varepsilon}$ and $\Nm_{a,\varepsilon} - \Nm_{0,\varepsilon}$, respectively, which we will do in what follows.

\section{Expansion of the external DtN operator}
\label{exter}
\setcounter{equation}{0}

\begin{lemma}
\label{l3a}
 For the unperturbed problem, operator $(\Np_{\bk,0}-\Nm_{0,0})$ has simple zero eigenvalue if $\bk$ is non-exceptional point. Otherwise, zero eigenvalue has multiplicity $n$ (where $n$ is the order of the exceptional point). The space  of eigenfunctions is spanned by
\begin{equation}\label{psi0}
\psi_{s}:=\E^{-\I (\bk - \mb_s) \cdot \br}|_{r=R}, \quad 1 \leq s \leq n,
\end{equation}
where $\mb_1 = \z$, points $\mb_s\in \mathbb Z^3_b$, $s > 1$, have been introduced in the definition of the exceptional point.
\end{lemma}
\begin{proof} Inclusionless problem \rf{Hz1},\rf{bc1a} has form \rf{unp}. Lemma \ref{lex} states that the solution space of the latter problem is spanned by functions \rf{psi0}. This solution space coincides with the kernel of the operator $(\Np_{\bk,0}-\Nm_{0,0})$ due to Lemma  \ref{l1a}. The proof is complete.
\end{proof}

We denote by $\mathscr E$ the finite-dimensional space spanned by functions $\psi_{s}$, $1\leq s \leq n$, and we denote by $\mathscr E_{\bot,1}, \mathscr E_{\bot,0}$ the subspaces in $\Ht(\de B_R)$ and $\Ho (\de B_R)$, respectively, that consist of the functions orthogonal in $L^2(\de B_R)$ to  $\mathscr E$.
We will write each element $\psi$ in the domain  $\Ht(\de B_R)$  and the range $\Ho(\de B_R)$ of operator \rf{npm} in the vector form $\psi=( \psi_{\mathscr E},\psi_\bot),$ where $\psi_{\mathscr E}$ is the projection in $L^{2}(\de B_R)$ of function $\psi$ into the space $\mathscr E$ ,
and $\psi_\bot$ is orthogonal to $\psi_{\mathscr E}$ in $L^{2}(\de B_R)$. Then, due to Lemmas \ref{l2}-\ref{l3a}, the unperturbed  operator \rf{npm} has the following matrix form:
\begin{equation}
\label{A}
 \Np_{\bk,0} - \Nm_{0,0} = \left(
 \begin{array}{cc}
  0 & 0 \\[2mm]
  0 & \Am
 \end{array}
 \right),
\end{equation}
where $\Am:\mathscr E_{\bot,1}\to \mathscr E_{\bot,0}$ is an isomorphism.

Since \rf{k-egv} depends smoothly on $\bk$ and $\varepsilon, ~|\varepsilon|<1$, operator $ \Np_{\bk,\varepsilon}$ is an infinitely smooth function of both arguments. Recall that $|\varepsilon|$ is small when $a$ is small, see Section \ref{out} . Thus the matrix representation of the operator $\Np_{\bk,\varepsilon} - \Nm_{0,\varepsilon}$ in the basis chosen for $\varepsilon = 0$ has the following form when $\varepsilon \to 0$:
\begin{equation}
\label{conc}
 \Np_{\bk,\varepsilon} - \Nm_{0,\varepsilon} = \left(
 \begin{array}{cc}
  \Cm\ei + O(\ei^2) & O(\ei) \\[2mm]
  O(\ei) &\Am + O(\ei)
 \end{array}
 \right)=\left(
 \begin{array}{cc}
  \Cm\ei + \ei^2 \Dm_{11} & \ei \Dm_{12} \\[2mm]
  \ei \Dm_{21} & \Am + \ei \Dm_{22}
 \end{array}
 \right), \quad \ei=(\kp -|\bk|)/|\bk|,
\end{equation}
where $\Cm$: $\mathscr E \to \mathscr E$ is a finite-dimensional operator, and operators $\Dm_{ij}=\Dm_{ij}(\ei,\bk)$ are infinitely smooth functions of the arguments. We fix the basis \rf{psi0} in $\mathscr E$ and identify operator $\Cm$ with its matrix representation in this basis. We evaluate $n\!\times\!n$ matrix $\Cm$ in the next Lemma.
\begin{lemma}\label{l33}
 Matrix $\Cm$ in the matrix expansion \rf{conc} of the operator
 $\Np_{\bk,\varepsilon} - \Nm_{0,\varepsilon}$ is equal to
 \begin{equation}
 \label{C}
  \Cm = 2|\bk|^2 |\Pi| \Im,
 \end{equation}
 where $\Im$ is the identity matrix.
\end{lemma}
\begin{remark}
 Hereafter the absolute value of a domain denotes its volume.
\end{remark}

\begin{proof}
If point $\bk$ is non-exceptional, then $\Cm$ is a number that is equal to the coefficient in the leading term of the asymptotics of
\begin{align}
C(\varepsilon) := \int_{\de B_R} \frac{\de u}{\de \n}\, \E^{\I \bk \cdot \x}\, \D S,  \quad   \frac{\de u}{\de \n} =(\Np_{\bk,\varepsilon}-\Nm_{0,\varepsilon})\,\E^{-\I \bk \cdot \x}.
 \label{Cc}
\end{align}
From \rf{Cc} it follows that $u=v-w$, where $v,w$ are solutions of the problems
 \begin{align}
 \label{v}
 \left( \Delta + k^2_{+} \right) v &= 0, \quad  \x\in \Pi \smallsetminus B_R,~~  \left\rrbracket  \E^{\I \bk \cdot \x} v(\x) \right\llbracket =0, ~~ \left. v\right|_{r=R}=\E^{-\I \bk \cdot \x},\\[2mm]
  \left( \Delta + k^2_{+} \right) w &= 0, \quad  \x\in  B_R,~~   \left. w\right|_{r=R}=\E^{-\I \bk \cdot \x}.
  \label{w1}
\end{align}
Let $v=\E^{-\I \bk \cdot \x}+\vt, ~w=\E^{-\I \bk \cdot \x}+\wt$. Then $\vt, ~\wt$ satisfy
 \begin{align}
 \left( \Delta + k^2_{+} \right) \vt &= \varepsilon_1 \E^{-\I \bk \cdot \x}, \quad  \x\in \Pi \smallsetminus B_R,~~  \left\rrbracket  \E^{\I \bk \cdot \x} \,\vt(\x) \right\llbracket =0, ~~ \left. \vt\right|_{r=R}=0,\\[2mm]
  \left( \Delta + k^2_{+} \right) \wt &= \varepsilon_1 \E^{-\I \bk \cdot \x}, \quad  \x\in  B_R,~~   \left. \wt\right|_{r=R}=0,
  \label{wt}
\end{align}
where $\varepsilon_1=|\bk|^2-k^2_{+} = -\varepsilon |\bk|(|\bk| + \kp) = -2\varepsilon |\bk|^2 + \O(\varepsilon^2)$. From here it follows that $\| \vt\|,~\| \wt\|=\O (\varepsilon)$, and therefore Green's second identity implies
\begin{align}
\label{gvt}
\int_{\de B_R} \frac{\de \vt}{\de \n}\, \E^{\I \bk \cdot \x}\, \D S&= -\int_{\Pi \smallsetminus B_R}  \left[\left( \Delta + |\bk|^2 \right) \vt\right]\E^{\I \bk \cdot \x}\, \D \x=-\varepsilon_1|\Pi \smallsetminus B_R|+\O(\varepsilon^2),
\\[2mm]
\int_{\de B_R} \frac{\de \wt}{\de \n}\, \E^{\I \bk \cdot \x}\, \D S&= \int_{B_R} \left[\left( \Delta + |\bk|^2 \right) \wt\right]\E^{\I \bk \cdot \x}\, \D \x=\varepsilon_1| B_R|+\O(\varepsilon^2).
\label{gwt}
\end{align}
Hence $C(\varepsilon) = 2\varepsilon |\bk|^2|\Pi|+\O(\varepsilon^2),$ and therefore $\Cm$ is given by \rf{C}.

Let now $\bk$ be exceptional, i.e. there are nonzero points $\bm m_j\in \mathbb Z^3_b$ such that $|\bk|=|\bk-\bm m_j|, ~ 2 \leq j \leq n$. Applying the same arguments we obtain that the diagonal entries of matrix $\Cm$ are equal to  $2|\bk|^2 |\Pi|$. For transparency, we show that the non-diagonal entries of $\Cm$ are zeros by considering the elements $C_{1,j}, \; j >1$. These elements are the coefficients for $\varepsilon$ in the asymptotics of the integral
\[
C_{1,j} (\varepsilon): = \int_{\de B_R} \frac{\de u}{\de \n}\, \E^{\I (\bk-\bm m_j) \cdot \x}\, \D S, \quad \mb_j \neq \z,
\]
where $\frac{\de u}{\de \n}$ is defined in \rf{Cc}. All the arguments above applied to evaluate $C(\varepsilon)$ can be repeated. Formulas \rf{v}-\rf{wt} remain valid but the integrands in \rf{gvt}-\rf{gwt} contain the exponents $\E^{\I (\bk-\bm m_j) \cdot \x},$ $ \mb_j \neq \z,$ instead of $\E^{\I \bk \cdot \x}$. This leads to
\[
C_{1,j} (\varepsilon)=\varepsilon\int_{\Pi}\E^{-\I \bm m_j \cdot \x}\, \D \x+\O(\varepsilon^2)=\O(\varepsilon^2).
\]
Therefore, $C_{1,j} = 0$.
\end{proof}

\section{Expansion of the inner DtN operator}
\label{inner}
\setcounter{equation}{0}

We consider the Dirichlet problem in the ball $B_R$ containing the inclusion $\Om$:
\begin{equation}
\left\{
\begin{array}{l}
\Delta u_{-} + k^2_{-}  u_{-} = 0, \quad u_{-} \in H^2 (\Om), \\[2mm]
\Delta u_{+} + k^2_{+}  u_{+} = 0, \quad u_{+}  \in H^2 (B_R \setminus {\Om}),
\end{array}
\right.
\label{eq1}
\end{equation}
\vspace{-3ex}
\begin{align}
 \left. u \right|_{ \de B_R} = \psi \in \Ht(\de B_R),
 \label{uR}
\end{align}
\vspace{-3ex}
\begin{align}
\label{bc1c}
\left\llbracket u  \right\rrbracket =0, \quad
\frac{\de u_{+} }{\de \n}  = \ep\,\frac{\de u_{-} }{\de \n},\quad  \x \in \de \Om, \quad \ep = \frac{\rho_{+}}{\rho_{-}}.
\end{align}
Notice that the second condition in \rf{bc1c} is equivalent to the corresponding jump condition in \rf{uin} or \rf{bc1a}.

It is more convenient for us to have inhomogeneity on $\de \Om$ rather than on the exterior boundary. Thus, we will be looking for the solution $u$ in the form
\begin{align}
u=u_0+\ut,
\label{ut}
\end{align}
where $u_0$ is the solution of the Dirichlet problem in the ball $B_R$ (without an inclusion):
\begin{align}
\Delta u_0+ k^2_{+} u_0 =0, \quad \x \in B_R ; \quad~~u_0|_{\de B_R}=\psi.
\label{u0}
\end{align}
Then $\ut$ satisfies
\begin{equation}
\left\{
\begin{array}{l}
\Delta \ut_{-} + k^2_{-}  \ut_{-} = (\kp^2 - \km^2) u_0, \quad \ut_{-} \in H^2 (\Om), \\[2mm]
\Delta \ut_{+} + k^2_{+} \ut_{+} = 0, \quad \ut_{+}  \in H^2 (B_R \setminus {\Om}),
\end{array}
\right.
\label{eq2a}
\end{equation}
\vspace{-2ex}
\begin{align}
 \left. \ut \right|_{ \de B_R} = 0,
 \label{uRa}
\end{align}
\vspace{-3ex}
\begin{align}
\label{bc1aa}
\ut_{-}   = \ut_{+}, \quad
\frac{\de \ut_{+} }{\de \n} - \ep\,\frac{\de \ut_{-} }{\de \n} =(\ep -1)\,\frac{\de u_0 }{\de \n} ,  \quad \x \in \de \Om.
\end{align}
The asymptotic expansion of the solution $u$ of \rf{eq1}-\rf{bc1c} as $a\to 0$ is rather complicated (it can be derived from \rf{sol} and Lemma \ref{l10}, see below). In fact, we will need only the asymptotics for the DtN map: $\Nm_{a,\varepsilon}: \left.\psi=u \right|_{\x\in B_R}\to \left.\frac{\de u }{\de \n}\right|_{\x\in B_R}$ which has the form of a simple power series. For the problem without inclusion, the image of this map $\Nm_{0,\varepsilon}$ is equal to $\left.\frac{\de u_0 }{\de \n}\right|_{\x\in B_R}$, and our goal is to find the difference
\begin{equation}\label{dn}
\Nm_{a,\varepsilon}-\Nm_{0,\varepsilon}:\psi\to \left.\frac{\de \ut }{\de \n}\right|_{\x\in B_R}.
\end{equation}
We also need the asymptotics of the quadratic form
\begin{align}
 \left( \left(\Nm_{a,\varepsilon}-\Nm_{0,\varepsilon}\right) \hpsi, \hpsi \right), \quad \hpsi \in \mathscr E,
 \label{M}
\end{align}
as $\varepsilon, a \to 0$. The matrix of this form in the basis $\hpsi_s,~ 1 \leq s \leq n,$ is denoted by $\M$.
We will use the same notation for the operator in $\mathscr E$ defined by the matrix $\M$.

Recall that $\Om = \Om (a)$ is obtained from an $a$-independent domain $\what{\Om}$ by the contraction with the coefficient $a^{-1}$, i.e. the transformation $\x \to a\si $ maps $\Om (a)$ into $\Omh$. We need a Taylor expansion of the function $u_0 (\x)$
\begin{align}
\label{u0_ser}
 u_0(\x) = u_0(\z) + {\bm c} \cdot \x + \oh\,{\Q \x}\cdot \x + \O \left( |\x|^3 \right), \quad |\x| \to 0, \quad {\bm c} =\nabla u_0 (\z), \quad \Q = \left\{ \frac{\de^2 u_0}{\de x_i \de x_j}\right\}_{\x = \z},
\end{align}
and the asymptotics at infinity of the auxiliary problem in the entire space $\R^3_\si$:
\begin{align}
 \label{V0}
 \Delta V &= 0, \ \quad \si \in \mathbb{R}^3\setminus\de\Omh, \\[2mm]
 \label{Vasy}
 V &\sim \O \left( \frac{1}{|\si|}\right), \quad |\si| \to \infty, \\[0mm]
 \label{V0asy}
 V_{+} &= V_{-},  \quad \si \in \de\Omh, \\[2mm]
 \frac{\de V_{+}}{\de \n} - \ep \frac{\de V_{-}}{\de \n} &=  (1-\ep) ({\bm c} \cdot \n), \quad \si \in \de\Omh.
 \label{bc_xi}
\end{align}
Here $\Delta=\Delta_\si$ denotes the Laplacian in the rescaled variables $\si$,  $\n$  denotes the external normal
vector to the surface $\Omh$, $V=V(\si)$, $V_-=V|_{\Omh} \in H^2 (\Omh)$, $V_+=V|_{{\R}^3\setminus\Omh} \in H^2_{\rm loc} ( {\R}^3\setminus\Omh)$, and $\ep$ is defined in \rf{bc1c}.

We will show that the problem is uniquely solvable and has the following asymptotics as $|\si| \to \infty$:
\begin{align}
 V = -\frac{1}{4\pi} \, \frac{{\bm p}\cdot \si}{|\si|^3} + \O \left( \frac{1}{|\si|^3}\right), \quad{\rm where} \quad  \m = \int_{\de \Omh} \si\,\walpha_0 (\si)\, \D S_{\si}, \quad \walpha_0=\left.\left(\frac{\de V_{+}}{\de \n} - \frac{\de V_{-}}{\de \n}\right)\right|_{\si \in \Omh}.
 \label{ppp}
\end{align}
Clearly, function $V$ and therefore vector $\m$ (the dipole moment of the surface $\de \Omh$) are proportional to $|{\bm c}|$. If $u_0(\x)=\E^{-\I \bk \cdot \x}$, then $|{\bm c}|=|\bk|$, and we introduce the polarizability vector $\bm \chi={\bm \chi}(\dotk),$ $\dotk =  \bk/|\bk|,$ by
\begin{align}
  {\bm \chi} =  \frac{\I \,\m}{|\bk| |\Omh|}.
 \label{chi}
\end{align}
If point $\bk$ is exceptional and $\m=\m_s$ corresponds to the external field $\E^{\I(\bk - {\bm m}_s) \cdot \x}$, then still $|{\bm c}|=|\bk - {\mb }_s|=|\bk|$  and the polarizability vector $\bm \chi=\bm \chi_s(\frac{\bk - {\mb }_s}{|\bk - {\mb }_s|})$ is defined by
\begin{align}
  {\bm \chi_s} =  \frac{\I \,\m_s}{|\bk| |\Omh|}.
 \label{chi1}
\end{align}
\begin{theorem}
\label{t2}
\item[(1)]
If $R \notin \{R_i\}$ then the operator-function $ \Nm_{a,\varepsilon}-\Nm_{0,\varepsilon}:\Ht(\de B_R)\to \Ho(\de B_R)$ is infinitely smooth in $a$ and $\ei$ in a neighborhood of $a=\ei=0$, and its Taylor expansion in $a$ starts with $a^3$:
\begin{equation}
\label{et2}
 \Nm_{a,\varepsilon}-\Nm_{0,\varepsilon} \sim \sum_{j=3}^\infty \P_j (\ei)\, a^j, \quad a\to 0,
\end{equation}
where $\P_j$ are bounded operators.
\item[(2)]
The entries $M_{i,j}$, $1\leq i,j \leq n$, of the matrix $\M$ have the form
\begin{align}
M_{i,j} = |\Pi| |\bk|^2 \left( 1 - \frac{\gi}{\ge} +{\bm \chi}_{i} \cdot\frac{\bk - {\bm m}_j}{|\bk - {\bm m}_j|}\right)f + \O \left(a^4 + a^3 |\varepsilon| \right), \quad a, \varepsilon \to 0,
 \label{Mat2}
\end{align}
where $f = |\Om|/|\Pi|$ is the volume fraction of the inclusions.
\end{theorem}
The proof of the theorem is given in Appendix A.

\section{Main results}
\label{main}
\setcounter{equation}{0}

The next theorem provides the dispersion relation when $\bk$ is fixed and it is non-exceptional.
\begin{theorem}
 \label{t1}
 If $\bk$ is a non-exceptional Bloch vector then the dispersion relation is given by
 \begin{align}
 \om^2 \rhop \ge = |\bk|^2 \left[1 + \left( 1 - \frac{\gi}{\ge} \right)f + ({\bm \chi} \cdot \dotk\,)  f  \right] + \O (a^4), \quad a \to 0.
 \label{dr1}
\end{align}
Solution of the problem \rf{Hz1},\rf{bc1a} has the form
\begin{align}
 u(\x) = C \E^{-\I \bk \cdot \x} \left(1 + \ut(\x)\right), \quad \| \ut \|_{C^\infty (\Pi \setminus B_R)} = \O(a),
 \quad a \to 0,
\end{align}
where $C$ is a constant.
\end{theorem}
Next we study problem \rf{Hz1},\rf{bc1a} when the Bloch vectors $\bk$ is fixed and it is exceptional of order $n$. Let $\M^0$ be the matrix with the entries
 \begin{align}
M^0_{i,j} =  1 - \frac{\gi}{\ge} + {\bm \chi}_{i} \cdot \frac{\bk - {\bm m}_j}{|\bk - {\bm m}_j|},
 \label{M0}
\end{align}
so that
\begin{align}
\label{M1}
\M = |\Pi| |\bk|^2 \M^0 f + \O \left(a^4 + a^3 |\varepsilon| \right).
\end{align}
 From Lemma \ref{l2} it follows that matrix $\M$ defined after \rf{M} is Hermitian. Then \rf{M1} implied that $\M^0$ is also Hermitian and, therefore, has real eigenvalues.
\begin{theorem}
 \label{t3}
 Let $\bk$ be an exceptional Bloch vector of order $n \geq 2$, $\mb_j\in\mathbb Z^3_b, 1
\leq j\leq n, |\bk-\mb_j|=|\bk|$. Assume that the eigenvalues $\lambda_s$, $1 \leq s \leq n,$ of matrix $\M^0$  are  distinct and  ${\bm \mu}_s = (\mu_{1,s},\ldots,\mu_{n,s})^\T$ are the corresponding eigenvectors.

 Consider a small neighborhood $I$ of the point $\om = |\bk| c_+ = |\bk|/\sqrt{\rhop \ge}$. Then, for small $a$, there are exactly $n$ values of the frequencies $\om = \om_s \in I$ for which the corresponding problem \rf{Hz1},\rf{bc1a} with $|k_\pm|^2 = \om_s^2 \rho_\pm \gamma_\pm$  has a non-trivial solution. The frequencies $\om_s$ satisfy the relation
 \begin{align}
 \om^2_s  \rhop \ge = |\bk|^2 \left(1 + \lambda_s f  \right) + \O (a^4), \quad a \to 0.
 \label{dr2}
\end{align}
The corresponding solutions of the problem \rf{Hz1},\rf{bc1a} have the form
\begin{align}\label{us1}
 u_s(\x) = C_s \left( \sum_{j=1}^n \mu_{j,s}\, \E^{-\I (\bk - \mb_j) \cdot \x} + \ut_{s}\right), \quad \| \ut_{s} \|_{C^\infty (\Pi \setminus B_R)} = \O(a),
 \quad a \to 0,
\end{align}
where $C_s$ are constants and
$\rrbracket  \E^{\I \bk \cdot \x}\, \ut_{s} \llbracket = 0$.
\end{theorem}
\begin{remark}
Thus, each solution $u_s (\x)$ is a cluster of waves propagating in different directions. If one takes a linear combination of solutions $u_s$ to obtain a wave propagating in one direction, then the combination would contain the terms with different time frequencies and these terms are solutions of different equations.
\end{remark}

Next, we consider Bloch vectors $\bk$ with a fixed direction. If the vector is not exceptional, then the dispersion relation is given in Theorem \ref{t2}. If $\bk$ coincides with an exceptional point on the chosen ray, then solution of the problem \rf{Hz1},\rf{bc1a} is studied in Theorem \ref{t3}. The following theorem concerns the case when $\bk$ may vary along the ray being close to an exceptional vector $\bk^\ast$ on that ray while the time frequency $\omega$ is fixed in such a way that $\kp =\omega/c_+=|\bk^\ast|$.
\begin{theorem}
\label{t4}
Let $\bk^\ast$ be an exceptional Bloch vector of order $n$ and $\om$ is chosen such that
$\om^2  \rhop \ge = |\bk^\ast|^2$. Assume that the eigenvalues $\lambda_s$, $1 \leq s \leq n,$ of matrix $\M^0$  are  distinct and  ${\bm \mu}_s = (\mu_{1,s},\ldots,\mu_{n,s})^\T$ are the corresponding eigenvectors.

Consider the interval $J$ of the wave vectors $\bk = (1+\delta)\bk^\ast$, $|\delta| \ll 1$.
Then, for small $a$, there are exactly $n$ values of the wave vectors ${\bk} = {\bk}_s \in J$ for which the problem \rf{Hz1},\rf{bc1a} with ${\bk} = {\bk}_s \in J$  has non-trivial solutions. The wave vectors ${\bk}_s$ satisfy the relation
 \begin{align}
 \bk_s = \bk^\ast \left(1 -\frac{1}{2} \lambda_s f  \right) + \O (a^4), \quad a \to 0.
 \label{dr3}
\end{align}
The corresponding solutions of the problem \rf{Hz1},\rf{bc1a} with $\kp = |\bk^\ast|$ and $\bk = \bk_s$ have the form
\begin{align}
\label{us}
 u_s(\x) = C_s \left( \sum_{j=1}^n \mu_{j,s}\, \E^{-\I (\bk_s - \mb_j) \cdot \x}  + \ut_{s}\right), \quad \| \ut_{s} \|_{C^\infty (\Pi \setminus B_R)} = \O(a),
 \quad a \to 0,
\end{align}
where $C_s$ are constants, $\mb_j\in\mathbb Z^3_b,~ 1
\leq j\leq n,~~ |\bk^\ast-\mb_j|=|\bk^\ast|$, and
$~\rrbracket  \E^{\I \bk_s \cdot \x}\, \ut_{s} \llbracket \,=\, \rrbracket  \E^{\I (\bk_s -\mb_j) \cdot \x}\, \ut_{s} \llbracket=0$.
\end{theorem}
\begin{remark} Solution \rf{us} is a cluster of waves propagating in different directions with the same spatial frequency $|\bk_s|$.

If one replaces $\bk_s$ in \rf{us} by $\bk^\ast$, then the estimate \rf{us} for $\ut_s$ remains valid and therefore one could construct a linear combination of functions $u_s$ to obtain a function proportional to
\[
v(\x)= \E^{-\I \bk^\ast \cdot \x}+\widetilde{v}(\x), \quad \| \widetilde{v} \|_{C^\infty (\Pi \setminus B_R)} = \O(a),
\]
which satisfies all the relations \rf{Hz1},\rf{bc1a} with $k_+=|\bk^\ast|$ except the Bloch periodicity condition. This function represents the wave propagating in one direction (of the vector $\bk^\ast$) but this wave is a combination of terms with different spatial frequencies $|\bk_s|$ and therefore the remainder term $\widetilde{v}$ is not small when the distance from the origin is large.
\end{remark}

\section{Proof of Theorems \ref{t1}-\ref{t4}}
\label{proof}
\setcounter{equation}{0}

\begin{proof}
We start with the proof of Theorem \ref{t3}.
Theorem \ref{t1} can be treated as a particular case of Theorem \ref{t3} with $n=0$. The same approach is used to prove Theorem \ref{t4}.

We use Lemma \ref{l1a} and reduce problem  \rf{Hz1},\rf{bc1a} to the equivalent  equation \rf{nn}. Due to Lemma \ref{l2}, the operator in the latter equation is Fredholm and symmetric. We rewrite   \rf{nn} in the form
\begin{equation}
\label{matr}
 (\Np_{\bk,\varepsilon}-\Nm_{0,\varepsilon})\psi - (\Nm_{a,\varepsilon}-\Nm_{0,\varepsilon})\psi=0.
\end{equation}
Operators in \rf{matr} act from $\Ht(\de B_R)$ to $\Ho (\de B_R)$. The dispersion relation is the relation between the parameters of the problem for which equation \rf{matr} has a non-trivial solution.

We split the domain and the range of operators into two orthogonal in $L^2(\de B_R)$ components where the first component $\mathscr E$ is a finite-dimensional space of functions spanned by functions $\psi_{s}$ defined in Lemma \ref{l3a}.
This allows us to rewrite \rf{matr} in a matrix form similar to that in \rf{A} and \rf{conc}.
In particular, if point $\bk$ is not exceptional then $\mathscr E$ is a one-dimensional space of functions proportional to $\E^{-\I \bk \cdot \x}$.
For any $\bk$, let $\psi = (\psi_{\mathscr E}, \psi_{\bot})^\T$ be the vector representation of function $\psi \in \Ht (\de B_R)$, i.e. $\psi_{\mathscr E}$  is the $L^2$-projection of $\psi$ into $\mathscr E$. From  \rf{conc} it follows that \rf{matr} has the form
\begin{equation}\label{matr1}
  \left(
 \begin{array}{cc}
  \Csbi\ei + \ei^2 \Dsbi_{11}-\B_{11}(\varepsilon,a) & \ei \Dsbi_{12}-\B_{12}(\varepsilon,a) \\[2mm]
  \ei \Dsbi_{21}-\B_{21}(\varepsilon,a) & \A + \ei \Dsbi_{22}-\B_{22}(\varepsilon,a)
 \end{array}
 \right)
  \left(
 \begin{array}{cc}
  \psi_{\mathscr E} \\[2mm]
  \psi_{\bot}
 \end{array}
 \right)=0, \quad \ei=(\kp -|\bk|)/|\bk|,
\end{equation}
where matrix elements $\B_{i,j}$ are defined by the operator $\Nm_{a,\varepsilon}-\Nm_{0,\varepsilon}$. Theorem \ref{t2} implies that $\|\B_{ij}\|=\O(a^3),~ \varepsilon,a\to 0.$ The element $\B_{1,1}$ coincides with the matrix $\M=\M(\varepsilon,a)$ defined after  \rf{M}.

{\it We will solve the fist equation in \rf{matr1} and find $\ei = \ei (a)$ for which this equation has a non-trivial solution. Then the second equation allows us to find $\kp = \kp (\ei)$ for which 
\rf{matr},\rf{matr1} and therefore \rf{Hz1},\rf{bc1a} have nontrivial solutions.}

Since operator $\A$ is invertible, the second equation of the system \rf{matr1} can be solved for $\phi_{\bot}$ yielding
\begin{align}
 \label{phi2_norm}
 \|\psi_{\bot}\| \leq c (|\ei| + a^3) \|\psi_{\mathscr E} \|.
\end{align}
This reduces \rf{matr},\rf{matr1} to an equation in the finite-dimensional space $\mathscr E$:
\begin{align*}
\label{kpe}
 [\Csbi\ei-\Msbi(\varepsilon,a)+O((|\ei|+a^3)^2)]\psi_{\mathscr E}=0.
\end{align*}
We substitute here formulas \rf{C} and \rf{M1} for $\Csbi$
and $\Msbi$, and obtain
\begin{equation}\label{drel}
\left[ 2\ei\Im-\M^0 f+\O\left(\ei^2+a^4\right)\right]\psi_{\mathscr E}=0.
\end{equation}
Since $\psi_{\mathscr E} = \sum_{j=1}^n \tau_j \psi_j$, last formula can be rewritten in terms of ${\bm \tau} = (\tau_1, \ldots, \tau_n)$
\begin{equation}\label{drel1}
\left[ 2\ei\Im-\M^0 f+\O\left(\ei^2+a^4\right)\right]{\bm \tau}=0.
\end{equation}
Here $\M^0$ is a matrix corresponding to the operator $\M^0$ in \rf{drel}, and its entries are given in \rf{M0}.
Since $\ep$ and $a$ are small, equation \rf{drel1} cannot have non-trivial solutions if $|\ei| \gg f$.
Matrix $\M^0$ is Hermitian and has different eigenvalues.
Hence, a non-trivial solution $\psi$ of \rf{matr} for small $\ei,a$ exists if and only if $2\ei/f$ approaches one of eigenvalues $\lambda_s$ of $\M^0$ as $a\to0$ and vector ${\bm \tau}$ approaches  the corresponding eigenvector ${\bm \mu}_s$ of $\M^0$. More precisely,
\begin{align}\label{2ei}
2\ei=\lambda_sf+\O(a^4), \quad {\bm \tau}={\bm \tau}_s=C_s{\bm \mu}_s (1+\O(a)), \quad a\to 0, \quad 1 \leq s \leq n.
\end{align}
For simplicity of formulas, we omit the subscript $s$ from $\varepsilon_s$.
Further, by multiplying the second relation in \rf{matr1} and the equality $\dst 2 = (\kp + |\bk|)/|\bk|  + \O(\ei)$, which follows from the second relation in \rf{matr1}, we obtain
\begin{align}
 2\ei =\frac{(k_{+} -|\bk|)(k_{+} +|\bk|)}{|\bk|^2} + \O(\varepsilon^2) = \frac{k_{+}^2 - |\bk|^2}{|\bk|^2}  + \O(\varepsilon^2) = \frac{k_{+}^2 - |\bk|^2}{|\bk|^2}  + \O(a^6).
 \label{Ce}
\end{align}
The last equality follows from \rf{2ei}.  Formulas \rf{2ei},\rf{Ce} imply that $\kp^2 = \om^2 \rhop \ge$ has $n$ distinct values for which $\om = \om_s,~ 1 \leq s \leq n,$ and
\begin{align}\label{77}
\om^2_s \rhop \ge - |\bk|^2=\lambda_s|\bk|^2f+\O(a^4), \quad
{\bm \tau}_s=C_s{\bm \mu}_s (1+\O(a)), \quad a \to 0.
\end{align}
The first relation above coincides with \rf{dr2}. To complete the proof of Theorem \ref{t3}, one needs only to justify \rf{us}.

Projection $\psi_{\mathscr E}$ of the solution $\psi$ of \rf{matr}, that corresponds to the vector ${\bm \tau}_s$, equals
\begin{align*}
\psi_{\mathscr E} = \psi_{\mathscr E,s} = C_s   \sum_{j=1}^n \left( \mu_{j,s}\, + \O(a) \right) \E^{-\I (\bk - \mb_j) \cdot \x}, \quad \x \in \de B_R.
\end{align*}
From hear and \rf{phi2_norm}, where an estimate of $\varepsilon$ can be taken from \rf{2ei}, it follows that solutions $\psi = \phi_s$, $1 \leq s \leq n$, of \rf{matr} have the form
\begin{align}
\label{78}
\phi_s = C_s \left(\sum_{j=1}^n \mu_{j,s}\, \E^{-\I (\bk - \mb_j) \cdot \x} + \tilde{\phi}_s \right), \quad \| \tilde{\phi}_s\|_{\Ht(\de B_R)} = \O(a), \quad a \to 0.
\end{align}
Function $u_s$ in the domain $\Pi \setminus B_R$ satisfies
\begin{align}
 \label{vv}
 \left( \Delta + k_{+}^2 \right) u_s &= 0, \quad  \x\in \Pi \smallsetminus B_R,~~  \left\rrbracket  \E^{\I \bk \cdot \x} u_s(\x) \right\llbracket =0, ~~ \left. u_s\right|_{r=R}=\phi_s,
 \end{align}
and $\ut_{s} $ is the solution of the problem
\begin{align}
 \label{vv1}
 \left( \Delta + k_{+}^2 \right) \ut_{s} &= (k_{+}^2-|\bk|^2) \sum_{j=1}^n \mu_{j,s}\, \E^{-\I (\bk - \mb_j) \cdot \x}, \quad  \x\in \Pi \smallsetminus B_R,~~  \left\rrbracket  \E^{\I \bk \cdot \x} \ut_{s}(\x) \right\llbracket =0, ~~ \left. \ut_{s}\right|_{r=R}=\tilde{\phi}_s .
 \end{align}

Denote by $\U_\ei:\Ht(\de B_R)\to H^{2}(\Pi \setminus B_R)$ the operator that maps each $h\in \Ht(\de B_R)$ into the solution of problem \rf{vv} with $\phi_s$ replaced by $h$, and denote by $\U_0$ similar operator when $ k_{+}^2$ in \rf{vv} is replaced by $|\bk|^2$. Due to the choice of $R\notin \{R_i\}$, operator $\U_0$ is bounded. Then operator $\U_\ei$ is also bounded when $a\ll 1$ since $k_{+}^2 -|\bk|^2=O(a^3)$ (see \rf{2ei} and the first two relations in \rf{Ce}). This implies the boundedness of a similar operator that maps a function on the right-hand side of equation \rf{vv1} into the solution of the problem with zero boundary conditions. Hence \rf{vv1} implies the estimate
\begin{align}
\label{hu2}
\| \ut_{s} \|_{H^2 (\Pi \setminus B_R)} = \O(a), \quad a \to 0,
\end{align}
 since the inhomogeneities in the equation and the boundary condition in \rf{vv1} have orders $O(a^3)$ and $O(a)$, respectively.
This estimate is valid for smaller values of $R$ and, therefore, the standard a priori estimates for elliptic equations in subdomains imply that the estimate \rf{hu2} is valid in $C^\infty$.
 Theorems \ref{t1}, \ref{t3} are proven.

All the arguments used to prove Theorems \ref{t1}, \ref{t3} up to formula \rf{2ei} are valid under conditions of Theorem \ref{t4}. Now, when $k_+$ is fixed and $|\bk|$ varies and is close to $k_+=|\bk^\ast|$, formula \rf{Ce} has the form
\begin{align}
 2\ei =\frac{\om^2 \rhop \ge - |\bk_s|^2}{|\bk_s|^2} + \O(\varepsilon^2), \quad 1 \leq s \leq n.
 \label{Ce1}
\end{align}
formula \rf{77} now becomes
\begin{align}
|\bk^\ast|^2- |\bk_s|^2=\lambda_s|\bk_s|^2f+\O(a^4), \quad
{\bm \tau}_s=C_s{\bm \mu}_s (1+\O(a)), \quad a \to 0.
\end{align}
The first relation can be solved for $|\bk_s|$:
\[
|\bk_s|=|\bk^\ast|\left(1-\oh \lambda_s f\right)+\O(a^4).
\]
Formula \rf{dr3} follows from the latter equality and the condition $\bk=\bk^\ast(1+\delta)$ assumed in the statement of Theorem \ref{t4}.
Formulas \rf{78}, \rf{vv} now take the form
\begin{align*}
\phi_s = C_s \left(\sum_{j=1}^n \mu_{j,s}\, \E^{-\I (\bk^\ast - \mb_j) \cdot \x} + \tilde{\phi}_s \right), \quad \| \tilde{\phi}_s\|_{\Ht(\de B_R)} = \O(a), \quad a \to 0, \quad 1 \leq s \leq n,
\end{align*}
\begin{align*}
 \left( \Delta + |\bk^\ast|^2 \right) u_s &= 0, \quad  \x\in \Pi \smallsetminus B_R,~~  \left\rrbracket  \E^{\I \bk_s \cdot \x} u_s(\x) \right\llbracket =0, ~~ \left. u_s\right|_{r=R}=\phi_s,
 \end{align*}
If $u_s$ is  presented in the form \rf{us} (to satisfy the Bloch periodicity condition and almost satisfy the boundary condition on $\de B_R$), then  $\ut_{s} $ is the solution of the problem
\begin{align*}
 \left( \Delta + |\bk^\ast|^2 \right) \ut_{s} &=  \sum_{j=1}^n (|\bk^\ast|^2-|\bk_s - \mb_j|^2) \,\mu_{j,s}\, \E^{-\I (\bk_s - \mb_j) \cdot \x}, \quad  \x\in \Pi \smallsetminus B_R, \\[2mm]
   \left\rrbracket  \E^{\I \bk \cdot \x} \ut_{s}(\x) \right\llbracket &=0, \quad \left. \ut_{s}\right|_{r=R}=\tilde{\phi}_s+\tilde{\phi}_s^0 ,
 \end{align*}
 where
 \[
 \tilde{\phi}_s^0=\sum_{j=1}^n \mu_{j,s}\, \E^{-\I (\bk^\ast - \mb_j) \cdot \x}-\sum_{j=1}^n \mu_{j,s}\, \E^{-\I (\bk_s - \mb_j) \cdot \x}=O(a^3), \quad a\to 0.
 \]
Since $|\bk_s - \bk^\ast| = \O(a^3)$ and $|\bk^\ast - \mb_j| = |\bk^\ast|$, the estimate \rf{us} for $\ut_{s}$ can be proved similarly to the corresponding estimate in the proof of Theorem \ref{t3}.
 \end{proof}

\section{Application to spherical inclusions}
\label{spher}
\setcounter{equation}{0}

We consider problem  \rf{Hz1}-\rf{bc1a} for a simple cubic lattice of spherical inclusions when $\Pi = [-\pi,\pi]^3$,   the reciprocal lattice is $\Z^3$, and $\Om$ is a ball of radius $a$.
\begin{description}[style=nextline]
\item[Dispersion at a non-exceptional point]
In the external field $u_0 (\x) = \E^{-\I \bk \cdot \x}$ with non-exceptional wave vector $\bk$, one can find from \rf{V0}-\rf{bc_xi} and \rf{alphah} that
\begin{align}
 \walpha_0 &= -3\I \, \frac{\ep -1}{\ep + 2}(\bk \cdot \n), \quad
 {\bm p} = -4\pi \I \bk \,\frac{\ep -1}{\ep + 2}, \quad
 \bchi = 3\dotk \,\frac{\ep -1}{\ep + 2}.
\end{align}
Then the dispersion relation \rf{dr1} becomes
\begin{align}
 \om^2  \ge = |\bk|^2 \nu_{+} \left[1 + \left( 1 - \frac{\gi}{\ge} \right)f + 3\frac{\nu_{-} - \nu_{+}}{\nu_{-} + 2\nu_{+}}  f  \right] + \O (a^4), 
\end{align}
where $\nu_{\pm} = 1/\rho_{\pm}$ is the specific volume. Using the relation $|\bk|^2 = \om^2 / c_{+}^2  + \O(f)$, the last formula can be written in the form
\begin{align}
\overline{\gamma} \omega^2 = \la \nu \ra |\bk|^2 + \O \left( a^4 \right),
 \label{k_sph}
\end{align}
where
\begin{align}
 \overline{\gamma} = \ge (1-f) + \gi f
\end{align}
is the average adiabatic bulk compressibility modulus of the medium and
the average specific volume $\la \nu \ra$ is determined by Maxwell's formula
\begin{align}
 \la \nu \ra = \nu_{+} \left(1 + 3\, \frac{\nu_{-} - \nu_{+}}{\nu_{-} + 2\nu_{+}}f \right).
 \label{}
\end{align}
Notice that the main term of the asymptotics \rf{k_sph} is continuous in $\bk$ and does not reveal spectral gaps when $\bk$ approaches an exceptional value (see also Conclusions).
\item[Dispersion at exceptional point of order two]
Consider now an exceptional Bloch vector $\bk^\ast = (\oh,\alpha,\beta)$ with $|\alpha|, |\beta| < \oh$. We will illustrate Theorem \ref{t4} and assume that $\om$ is fixed in such a way that $\kp = \frac{\om}{c_{+}}= |\bk^\ast|$ in \rf{Hz1},\rf{bc1a}. There is only one non-trivial integer-valued vector $\mb = (1,0,0)$ such that $|\bk^\ast| = |\bk^\ast - \mb|$, i.e. the order of $\bk^\ast$ is two. {\it In the absence of an inclusion, problem \rf{Hz1},\rf{bc1a} with $\kp = |\bk^\ast|$  and $\bk = \bk^\ast$ has two linearly independent solutions $\E^{-\I \bk^\ast \cdot \x}$ and $\E^{-\I (\bk^\ast - \mb) \cdot \x}$.} Evaluation of matrix $\M^0$ in \rf{M0} gives
\begin{align}
 \label{}
 \M^0 &= \left(
 \begin{array}{cc}
 1 - \dfrac{\gi}{\ge} + 3\,\dfrac{\ep -1}{\ep + 2} & 1 - \dfrac{\gi}{\ge} + 3\,\dfrac{\ep -1}{\ep + 2}\, \dfrac{-1 + 4\alpha^2 + 4\beta^2}{1 + 4\alpha^2 + 4\beta^2}\\[2mm]
 1 - \dfrac{\gi}{\ge} + 3\,\dfrac{\ep -1}{\ep + 2}\, \dfrac{-1 + 4\alpha^2 + 4\beta^2}{1 + 4\alpha^2 + 4\beta^2}& 1 - \dfrac{\gi}{\ge} + 3\,\dfrac{\ep -1}{\ep + 2}
 \end{array}
 \right).
\end{align}
The eigenvalues and eigenvectors of $\M^0$ are
\begin{alignat}{5}
&\lambda_1 &&= 6\,\dfrac{\ep -1}{\ep + 2} \,\frac{1}{1 + 4\alpha^2 + 4\beta^2}, \quad &&{\bm \mu}_1 &&= (-1,1)^\T \\[2mm]
&\lambda_2 &&= 2\left(1 - \dfrac{\gi}{\ge}\right) + 24\,\dfrac{\ep -1}{\ep + 2} \,\frac{\alpha^2 + \beta^2}{1 + 4\alpha^2 + 4\beta^2}, \quad &&{\bm \mu}_2 &&= (1,1)^\T.
 \label{}
\end{alignat}
From Theorem \ref{t4} it follows that the Bloch vectors of the perturbed medium have the form
\begin{align}
\label{bk0}
\bk_1 &= \bk^\ast \left(1 - \oh\, \lambda_1 f \right) + \O(a^4), \\[2mm]
\bk_2 &= \bk^\ast \left(1 - \oh\, \lambda_2 f \right) + \O(a^4),
\label{bk1}
\end{align}
{\it There are two linearly independent solutions $u_1$ and $u_2$ of the perturbed problem when $\om = |\bk^\ast| c_{+}$. Each solution is a cluster of two waves} that have the following asymptotics outside of a neighborhood of the inclusions
\begin{align}
\label{v0}
u_1 (\x) &= - \E^{-\I \bk_1 \cdot \x} + \E^{-\I (\bk_1 - \mb) \cdot \x} + \O (a^3), \\[2mm]
u_2 (\x) &= \E^{-\I \bk_2 \cdot \x} + \E^{-\I (\bk_2 - \mb) \cdot \x} + \O (a^3),
\label{v1}
\end{align}
where the estimates of the remainder terms do not depend on the cell of periodicity. Cluster $u_1$ is a combination of two waves propagating in different directions with the same spatial frequency $|\bk_1| = |\bk_1 - \mb|$ that can be found from \rf{bk0}. Spatial frequencies of the waves in cluster $u_2$ are also equal but different from those in $u_1$.

One could consider function $v = u_2 - u_1$ in a hope that the terms with vector $\mb$ would be canceled with the accuracy of $\O(a^3)$, so that
\begin{align}
 v = \E^{-\I \bk_2 \cdot \x} + \E^{-\I \bk_1 \cdot \x} +\O(a^3)
 \label{vv2}
\end{align}
would represent a wave propagating in the direction of $\bk^\ast$. However,
the exponents in \rf{v0}, \rf{v1} have different spatial frequencies, and the difference of the exponents containing vector $\mb$ is small only when $|\x|$ is bounded. Hence, \rf{vv2} is not valid when $|\x|$ is large.

Choosing another vector from the same exceptional set does not change the set of solutions $\{u_1, u_2\}$. If one starts with the exceptional Bloch vector $\tilde{\bk}^\ast = (-\oh,\alpha,\beta)$, then $\mb = (-1,0,0)$, matrix $\M^0$ and its eigenvalues do not change, and the Bloch vectors in the perturbed medium become
\begin{align}
\label{ut0}
\tilde{\bk}_1 &= \tilde{\bk}^\ast \left(1 - \oh\, \lambda_1 f \right) + \O(a^4), \\[2mm]
\tilde{\bk}_2 &= \tilde{\bk}^\ast \left(1 - \oh\, \lambda_2 f \right) + \O(a^4).
\label{ut1}
\end{align}
The corresponding solutions outside of the neighborhood of the inclusions are
\begin{align}
\tilde{u}_1 (\x) &= - \E^{-\I \tilde{\bk}_1 \cdot \x} + \E^{-\I (\tilde{\bk}_1 - \tilde{\mb}) \cdot \x} + \O (a^3), \\[2mm]
\tilde{u}_2 (\x) &= \E^{-\I \tilde{\bk}_2 \cdot \x} + \E^{-\I (\tilde{\bk}_2 - \tilde{\mb}) \cdot \x} + \O (a^3),
\end{align}
and $\tilde{u}_1 (\x) = -u_1 (\x)$, $\tilde{u}_2 (\x) = u_2 (\x)$.
\item[Dispersion at exceptional point of order four]
The exceptional Bloch vector ${\bk}^\ast = (\oh,\frac{1}{3},\frac{2}{3})$ has order four with the other three exceptional vectors
$(-\oh,\frac{1}{3},\frac{2}{3})$, ($\mb_2 = (1,0,0)$),  $(\oh,-\frac{2}{3},-\frac{1}{3})$, ($\mb_3 = (0,1,1)$, $(-\oh,-\frac{2}{3},-\frac{1}{3})$, ($\mb_4 = (1,1,1)$). Calculation of the eigenvalues and eigenvectors of the matrix $\M^0$ gives
\begin{alignat}{5}
&\lambda_1 &&= 0,
\quad &&{\bm \mu}_1&& = (1,-1,-1,1)^\T, \\[2mm]
&\lambda_2 &&= \frac{108}{29}\,\frac{\ep -1}{\ep + 2},
\quad &&{\bm \mu}_2&& = (-1,1,-1,1)^\T, \\[2mm]
&\lambda_3 &&= \frac{216}{29}\,\frac{\ep -1}{\ep + 2},
\quad &&{\bm \mu}_3&& = (-1,-1,1,1)^\T, \\[2mm]
&\lambda_4 &&= 4\left(1 - \frac{\ge}{\gi}\right)+\frac{24}{29}\,\frac{\ep -1}{\ep + 2},
\quad &&{\bm \mu}_4&& = (1,1,1,1)^\T.
\end{alignat}
Then Theorem \ref{t4} implies that the perturbed Bloch vectors have the form
\begin{align}
\bk_s &= \bk^\ast \left(1 - \oh\, \lambda_s f \right) + \O(a^4), \quad 1 \leq s \leq 4.
\end{align}
By Theorem \ref{t4} there are four linearly independent solutions of \rf{Hz1},\rf{bc1a},  each of which is a cluster consisting of four waves
\begin{align}
\label{}
u_1 (\x) &= \E^{-\I \bk_1 \cdot \x} - \E^{-\I (\bk_1 - \mb_2) \cdot \x} -  \E^{-\I (\bk_1 - \mb_3) \cdot \x} + \E^{-\I (\bk_1 - \mb_4) \cdot \x} + \O (a^3), \\[2mm]
u_2 (\x) &= -\E^{-\I \bk_2 \cdot \x} + \E^{-\I (\bk_2 - \mb_2) \cdot \x} -  \E^{-\I (\bk_2 - \mb_3) \cdot \x} + \E^{-\I (\bk_2 - \mb_4) \cdot \x} + \O (a^3), \\[2mm]
u_3 (\x) &= -\E^{-\I \bk_3 \cdot \x} - \E^{-\I (\bk_3 - \mb_2) \cdot \x} +  \E^{-\I (\bk_3 - \mb_3) \cdot \x} + \E^{-\I (\bk_3 - \mb_4) \cdot \x} + \O (a^3),\\[2mm]
u_4 (\x) &= \E^{-\I \bk_4 \cdot \x} + \E^{-\I (\bk_4 - \mb_2) \cdot \x} +  \E^{-\I (\bk_4 - \mb_3) \cdot \x} + \E^{-\I (\bk_4 - \mb_4) \cdot \x} + \O (a^3).
\label{}
\end{align}

\end{description}

\section{Conclusions}
\label{concl}

We derived and rigorously justify asymptotic expansions for Bloch waves in periodic media with small inclusions of arbitrary shape and transmission boundary conditions. When the wave vector $\bk$ is not exceptional the wave propagates in a certain direction. The rigorous approach reveals that there are exceptional wave vectors  $\bk$ for which solution space is finite-dimensional, and each solution of the problem is a cluster of Bloch waves propagating in different directions or one direction with different frequencies.

Asymptotics of the solution provides the dispersion relation when $\bk$ is not exceptional. The dispersion relation cannot be defined uniquely when the solution is a cluster of Bloch waves.

Asymptotics as $a \to 0$ cannot be applied immediately for an investigation of spectral gaps, and it is not only due to the existence of clusters. For example, the existence of two different frequencies on the edge of the Brillouin zone does not guarantee that the surfaces showing the spectrum of the problem as a function of the wave vector are separated. The difference of the values could be explained by an intersection of these surfaces inside the Brillouin zone. To justify the existence of a gap at an exceptional point $\bk^\ast$, one needs to find the spectrum in a neighborhood of $\bk^\ast$ and pass to the limit when $\bk \to \bk^\ast$. Unfortunately, asymptotic expansion of the solution in $a$ is not uniform in $\bk$, and the limit as $\bk \to \bk^\ast$ of the main term of asymptotics in $a$ does not produce a spectral gap.

\section*{Appendix A.}

Below we provide the proof of Theorem \ref{t2} which is split into several sections.

\subsection*{Reduction of the problem to a system of integral equations}
\setcounter{equation}{0}
\setcounter{section}{1}
\renewcommand{\theequation}{\Alph{section}.\arabic{equation}}

We reduce \rf{eq2a}-\rf{bc1aa} to the integral equations using the potentials with the kernel defined by Green's function of the Dirichlet problem for the Helmholtz equation.
We use the following notation.
\begin{itemize}
 \item
$G(\x, \y)$ denotes Green's function of the Dirichlet problem in the ball $B_R$ without the inclusion:
\begin{align}\label{kern}
 \left( \Delta + k^2_{+} \right) G &= \delta (\x - \y), ~~\quad
 \left. G \right|_{\de B_R} = 0.
\end{align}
 \item
$\G$ denotes the operator with the kernel $G(\x,\y)$ defined on functions in $\Om$ (For operators, we often will use a bold version of the same letter which is used for its integral kernel)
\begin{align}
 \left( \G w\right) (\x) = \int_{\Om} G(\x,\y)\, w(\y)\,\D \y, \quad \x \in B_R.
 \label{Gop}
\end{align}
 \item
 We denote by $u^\prime_{\pm}$ the limiting values on $\de \Om$ of the normal derivatives $\frac{\de u}{\de \n}$ when $\x \to \de \Om$ from outside and inside of $\de \Om$, respectively. We omit $\pm$ if the normal derivative is continuous.
 \item
 $\G^\prime v = (\G \,v)^\prime$.
 \item
 $\G_{\de} \alpha$ denotes the surface potential
 \begin{align}
  \left( \G_{\de} \alpha \right) (\x) = \int_{\de \Om} G(\x,\y)\, \alpha (\y)\, \D S_{\y}, \quad \x \in B_R.
  \label{Gd}
 \end{align}
 \item
 $\dst \left(\G_{\de}\right)_{\pm}^\prime \alpha = \dst \left(\G_{\de} \alpha \right)_{\pm}^\prime$ are the limiting values of the normal derivatives on $\de \Om$ from outside and inside of $\de \Om$, respectively.
\end{itemize}

We will be looking for a solution of \rf{eq2a}-\rf{bc1aa} in the form
\begin{align}
 \ut &= \int_{\Om} G(\x, \y) \, w(\y)\, \D {\y} +\int_{\de \Om} G(\x, \y) \, \alpha(\y)\, \D S_{\y} = \G w  +  \G_{\de} \alpha,
 \quad ~ w \in L^2(\Om ),~\alpha \in \Ho (\de \Om ).
\label{sol}
\end{align}

\begin{lemma}
\label{lemma1}
 Formula \rf{sol} provides a one-to-one correspondence between solutions $u$ of \rf{eq2a}-\rf{bc1aa} and solutions of the system
 \begin{align}
 \label{w}
 w + q\G w + q\G_{\de} \alpha &= -qu_0, \quad \x \in \Om,  \\[2mm]
 \oh\, \alpha  - \varkappa \T \alpha - \varkappa \G^\prime w &= \varkappa u_0^\prime, \quad \x \in \de \Om, \quad ~~  w \in L^2(\Om ),~\alpha \in \Ho (\de \Om ),
 \label{alpha}
\end{align}
where $q = \km^2 - \kp^2, ~\varkappa = \frac{\rho_{+} - \rho_{-}}{\rho_{+} + \rho_{-}} = \frac{\ep -1}{\ep + 1}$.
\end{lemma}

\begin{proof}
Let $w, \alpha$ be a solution of \rf{w}-\rf{alpha}. We need to show that $\ut$ defined by \rf{sol} is a solutions of \rf{eq2a}-\rf{bc1aa}. Substitution of \rf{sol} into \rf{eq2a} leads to \rf{w}. Thus, \rf{eq2a} is valid. Further, \rf{sol} satisfies \rf{uRa}.

To rewrite \rf{bc1aa} in terms of $w$ and $\alpha$, we single out the singularity of $G(\x, \y)$:
\begin{align}\label{kern1}
 G(\x, \y) = - \frac{\cos k_{+} |\x - \y|}{4\pi |\x - \y|} + \tilde{G}, \quad \tilde{G} \in C^\infty, \quad |\y|\ll1, \quad \x\in  B_R,
\end{align}
and expand the first term in a power series in $|\x - \y|$:
\begin{align}
 G(\x, \y) \sim  - \frac{1}{4\pi |\x - \y|} + \sum_{n=0}^\infty a_n |\x - \y|^{2n+1} + \tilde{G}.
 \label{G}
\end{align}
From \rf{G} it follows that $\G_{\de} \alpha,~\G w, ~\G^\prime \alpha$ have the same properties as the standard simple layer, volume, and double layer potentials, respectively, in the potential theory. In particular,
\begin{align}
 \G^\prime_{\de \pm} \alpha = \pm \frac{\alpha}{2} + \T \alpha,
 \label{Gp}
\end{align}
where $\T$ is an integral operator with the kernel $T(\x, \y) = \dfrac{\de G}{\de \n}$, i.e.
\begin{align}
  \left( \T  \alpha \right) (\x) = \int_{\de \Om} \frac{\de G(\x,\y)}{\de \n}\, \alpha (\y)\, \D S_{\y},
  \label{T}
 \end{align}
\begin{align}
 T(\x,\y) \sim \frac{(\x-\y) \cdot \n_{\x}}{4\pi |\x -\y|^3} + \sum_{n=0}^\infty (2n+1) a_n \left[(\x-\y) \cdot \n_{\x}\right] |\x - \y|^{2n-1} + \tilde{T}, \quad \tilde{T} \in C^\infty.
 \label{T1}
\end{align}
Hence from \rf{sol} and \rf{Gp} we have
\begin{align}
 \label{upp}
 \ut^\prime_{+} &= \G^\prime w+ \frac{\alpha}{2} + \T \alpha, \\[2mm]
 \ut^\prime_{-} &= \G^\prime w- \frac{\alpha}{2} + \T \alpha.
 \label{upm}
\end{align}
Thus function \rf{sol} satisfies the first relation in \rf{bc1aa}, and the condition of the jump of the normal derivative in \rf{bc1aa} is equivalent to \rf{alpha}.
Hence, function $\ut$ satisfies \rf{eq2a}-\rf{bc1aa}.
Obviously, \rf{sol} defines $\ut$ by $w,\alpha$ uniquely.

Conversely, let $\ut$ be a solution of \rf{eq2a}-\rf{bc1aa}. Then
\begin{equation}\label{uw}
\Delta \ut + \kp^2  \ut = w+\alpha\delta(\partial\Omega),\quad \x \in B_R, \quad
w = \left\{
\begin{array}{cl}
 (\kp^2 - \km^2)(\ut+u_0), & \x \in \Om, \\[2mm]
 0, &\x \notin \Om,
\end{array}
\right.
\end{equation}
where  $\delta(\de \Om)$ is the delta-function on the surface $\de \Om$, and the coefficient $\alpha$ is equal to the jump of the normal derivative of $\ut$ on $\de \Om$. Hence, $\ut$ can be represented as $ \G (w + \alpha\delta(\partial\Omega))$, i.e.  by \rf{sol}. It was shown that $w,\alpha$ satisfy \rf{w}-\rf{alpha} in this case. It is also true that \rf{sol} defines $w,\alpha$ by $\ut$ uniquely. Indeed, if $\ut$ has form \rf{sol} and $\ut=0$, then the jump $\alpha$ of the normal derivative of $\ut$ on $\de\Om$ is zero, and the application of operator $\Delta+\kp^2$ to  \rf{sol} implies that $w=0$.
\end{proof}

\subsection*{Integral equations in new variables}

We introduce new variables in which $\Om = \Om (a)$ becomes an $a$-independent domain $\what{\Om}$ and introduce some notation related to this change:
\begin{align}
 \x &\to \si, \quad \si = \x/a,\\[2mm]
\ut(\x) &\to v(\si):=\ut(a\si), \\[2mm]
 \Om \in \mathbb{R}^3_\x  &\to \Omh \in\mathbb{R}^3_{\si}, \\[2mm]
 B_R &\to B_{R/a}, \\
 q(\x) &\to
 \qh (\si) = \left\{
 \begin{array}{cc}
  \km^2 - \kp^2, & \si \in \Omh, \\[2mm]
  0, & \si \notin \Omh.
 \end{array}
 \right.
 \label{qh}
\end{align}
We use $v$ to denote function $\ut$ after rescaling $\si = \x/a$,
and use the ``hat'' to denote other functions and operators in the new variables. In variables $\si$ system \rf{w}-\rf{alpha} becomes
\begin{equation}
 \label{sys1}
\Am \left(
    \begin{array}{c}
      \wh \\
      \walpha \\
    \end{array}
  \right)=\left(
            \begin{array}{c}
              -q\uh_0 \\
              \varkappa\uh_0' \\
            \end{array}
          \right), \quad {\rm where}  \quad
          \Am =\left(
              \begin{array}{cc}
               \Id + q \what{\G}  &q\what{\G}_{\de} \\
                -\varkappa \what{\G}^\prime &\dfrac{1}{2} \Id-\varkappa \what{\T} \\
              \end{array}
            \right),
\end{equation}
 and the matrix operator $\Am$ acts in $L^2(\Omh)\times \Ho(\de\Omh)$.

System \rf{sys1}  is defined on functions whose domains of definition are $a-$independent. Moreover, formula \rf{G}, where function $\tilde{G}$ is infinitely smooth,  and analyticity of $u_0(\x)$ in $|\x|$ for small $|\x|$ immediately imply the validity of the following statement.
\begin{lemma}
\label{l6}
Operator $\Am$ and its components
\[
\what{\G}:L^2(\Omh)\to L^2(\Omh),~~~
\what{\G}_{\de}:\Ho(\de\Omh)\to L^2(\Omh),~~~
\what{\G}^\prime:L^2(\Omh)\to \Ho(\de\Omh),~~~
\what{\T}:\Ho(\de\Omh)\to \Ho(\de\Omh),
\]
as well as the right-hand sides in \rf{sys1} are analytic in $a$ and $\varepsilon$ in a neighborhood of the point $a = \varepsilon=0$.
\end{lemma}
\begin{remark}
The main terms of the Taylor expansions in $a$ for $\Am$ and the right-hand side in \rf{sys1} are particularly important, and they are derived below.
\end{remark}
\begin{proof}
Operators and functions mentioned in the lemma are defined through solutions of the Dirichlet problems in the ball $B_R$ for the Helmholtz  operator $\Delta+k_+$ where $k_+=|\bk|(1+\varepsilon)$ and $\bk$ is fixed. Hence, the analyticity in $\varepsilon$ follows from the analyticity of the coefficient of the equation. Let us show the analyticity in $a$. We start with the operator $\what{\G}$. The integral kernel $G(x,y)$ of operator $\G$ is defined by \rf{kern}, \rf{kern1}. Hence $\tilde{G}$ is the solution to the problem
\begin{align*}
 \left( \Delta + k^2_{+} \right) \tilde{G} &= 0, ~~\quad
 \left. \tilde{G} \right|_{\de B_R} = \frac{\cos k_{+} |\x - \y|}{4\pi |\x - \y|} .
\end{align*}
The boundary function here is analytic in $y$ when $|y|<R$, and therefore the solution $\tilde{G}$ is analytic in $y$. It is also analytic in $|x|$ when $|x|<R$ since this is true for an arbitrary solution of the Helmholtz equation in the ball $B_R$. Hence $\tilde{G}(a\si,a\et)$ is analytic in $a$ when $\si,\et\in \Omh, a\ll 1$. Thus \rf{G} implies that operator $\what{\G}$ with the integral kernel $\what{G}(\si,\et)=a^3 G(a\si,a\et)$ is analytic in $a$ and has the following Taylor expansion at $a=0$:
\begin{align*}
 \what{\G} = a^2 \sum_{n=0}^\infty \what{\G}_n \,a^n, \quad \what{G}_0 (\si,\et) = -\frac{1}{4\pi |\si - \et|},
 \quad \si, \et \in \Omh.
\end{align*}

The same approach works for all other operators mentioned in Lemma \ref{l6}. In particular, since $\what{G}_{\de}(\si,\et)=a^2 G_{\de}(a\si,a\et)$, we have
\begin{align*}
\what{\G}_{\de} = a {\what {\F}}_{\de,1} + \sum_{n=3}^\infty {\what{\F}}_{\de,n} \,a^n, \quad {\what F}_{\de,1} (\si,\et) = -\frac{1}{4\pi |\si - \et|}, \quad \quad \si \in \Omh, \quad \et \in \de \Omh.
\end{align*}
Similarly,

\begin{align}
 \what{\T} = \what{\T}_0 + \sum_{n=2}^\infty a^n \what{\T}_n,  ~\quad \what{T}_0 (\si,\et) = \frac{(\si-\et) \cdot \n_{\si}}{4\pi |\si -\et|^3}, \quad \si, \et \in \de \Omh,
 \label{T0h}
\end{align}
\begin{align}
\label{Gpp}
 \what{\G}^\prime = a\what{\G}^\prime_1+ \sum_{n=3}^\infty \what{\G}^\prime_n a^n, \quad
 \what{\G}^\prime_1 \wh =  \int_{\what{\Om}} \what{T}_0(\si,\et)\, \wh(\et)\, \D \et, \quad \si \in \de \what{\Om}.
\end{align}
Analyticity in $a$ of the right-hand side in \rf{sys1}  is a consequence of the analyticity in $|\x|$ of the solution $u_0$ of the Helmholtz equation in the ball $B_R$ (see \rf{u0}).
\end{proof}

Using Lemma \ref{l6} and the expansions obtained in its proof, we can rewrite system \rf{sys1} in the form of power series:
\begin{equation}\label{sysC}
\sum_{n=0}^\infty a^nA_n\left(
                          \begin{array}{c}
                            \wh \\
                           \walpha\\
                          \end{array}
                        \right)=\sum_{n=0}^\infty a^n\left(
                          \begin{array}{c}
                           g_n\\
                           h_n\\
                          \end{array}
                        \right),  \quad ~~ A_0=
\left(
  \begin{array}{cc}
  \Id & 0 \\
    0 & \frac{1}{2}\Id - \varkappa \what{\T}_0 \\
  \end{array}
\right), \quad ~~  A_1=
\left(
  \begin{array}{cc}
    0 & q{\what {\F}}_{\de,1} \\
  -\varkappa \what{\G}^\prime_1 &0 \\
  \end{array}
\right),
\end{equation}
where $g_n,~h_n$ are the coefficients in the Taylor expansions of the functions on the right-hand side of \rf{sys1} :
\begin{align}
\label{u0exp}
 -qu_0(a\si) &= \sum_{n=0}^\infty a^ng_n(\si), \quad g_0 =-q\, u_0 (\z),
 \quad g_1 =-q\, {\bm c} \cdot \si, \ldots,\\[2mm]
 \varkappa u_0^\prime (a\si) &= \sum_{n=0}^\infty a^nh_n(\si), \quad
 h_0= \varkappa \,\n \cdot {\bm c}, \quad
 h_1= \varkappa \,\Q \si \cdot \n ,\ldots.
 \label{u1exp}
\end{align}
Here ${\bm c}, \Q$ are defined by the Taylor expansion \rf{u0_ser} of function $u_0$.
It will be shown below (see Lemma \ref{lemma3}) that the solution $\wh, \walpha$ can also be represented as a power series.
Then the main terms of the solution of \rf{sysC} satisfy
\begin{align}
\left(
\begin{array}{cc}
 \Id & 0 \\[2mm]
 0 & \frac{1}{2}\,\Id - \varkappa \what{\T}_0\
\end{array}
\right)
\left(
\begin{array}{c}
 \wh_0 \\[2mm]
 \walpha_0
\end{array}
\right) = \left(
\begin{array}{c}
 g_0 \\[2mm]
 h_0
\end{array}
\right),
 \label{approx0}
\end{align}
and
\begin{align}
\left(
\begin{array}{cc}
 \Id & 0 \\[2mm]
 0 & \frac{1}{2}\,\Id - \varkappa \what{\T}_0\
\end{array}
\right)
\left(
\begin{array}{c}
 \wh_1 \\[2mm]
 \walpha_1
\end{array}
\right) = \left(
\begin{array}{c}
 g_1 - q{\what {\F}}_{\de,1}\\[2mm]
 h_1 + \varkappa \what{\G}^\prime_1 \wh_0
\end{array}
\right).
 \label{approx1}
\end{align}
In particular,
\begin{align}
 \wh_0 = -q\, u_0 (\z).
 \label{w0}
\end{align}

\subsection*{Asymptotics of solutions of integral equations}

 The following statement allows one to find solutions of \rf{sysC}. Consider the elliptic problem \rf{V0}-\rf{V0asy} in the whole space with a general inhomogeneity in the last equation:
 \begin{align}
 \left(\frac{\de V_+}{\de \n} - \ep \,\frac{\de V_{-}}{\de \n}\right)\left|_{\si \in \de\Omh}\right. &= (1+\ep)f \in \Ho(\de\Omh).
 \label{sys72}
\end{align}
\begin{lemma}\label{l3} Problem \rf{V0}-\rf{V0asy}, \rf{sys72} is uniquely solvable. Operator $ \frac{1}{2}\Id - \varkappa \what{\T}_0 :\Ho(\de\Omh)\to \Ho(\de\Omh)$ is invertible, and the
  solution of the equation
  \begin{equation}\label{ttt}
\left(\oh\,\Id - \varkappa \what{\T}_0\right)\walpha=f\in \Ho(\de\Omh),
  \end{equation}
  is equal to
\begin{align}
 \walpha = \left. \left( \frac{\de V_+}{\de \n} - \frac{\de V_{-}}{\de \n} \right)\right|_{\si \in \de\Omh}.
 \label{alphah}
\end{align}
\label{lemma2}
\end{lemma}

\begin{proof}
We are looking for a solution of \rf{V0}-\rf{V0asy}, \rf{sys72} in the form of a simple layer  potential
 \begin{align}
  V&= -\int_{\de \Omh}\frac{1}{4\pi |\si - \et|}\, \walpha (\et)\, \D S_{\et}, \quad \walpha\in \Ho (\de \Omh).
  \label{vh}
 \end{align}
From the standard potential theory it follows that \rf{vh} satisfies \rf{V0}-\rf{V0asy}, \rf{sys72} if and only if \rf{ttt} holds. The potential theory also implies \rf{alphah}. Since solutions $V$ of \rf{V0}-\rf{V0asy}, \rf{sys72} satisfy $\Delta V=\walpha\delta(\de\Omh)$, where $\delta(\de\Omh)$ is the delta function on the surface $\de\Omh$ and $\walpha$ is the jump of the normal derivative of $V$ on that surface, each solution $V$ can be represented in the form \rf{vh}.  Hence, \rf{vh} establishes the one-to-one correspondence between solutions of \rf{V0}-\rf{V0asy}, \rf{sys72} and solutions of \rf{ttt}. Since both problems are Fredholm, the proof will be completed if the uniqueness of the solution to \rf{V0}-\rf{V0asy}, \rf{sys72} is shown.

From \rf{vh} it follows that
\[
|V|=O(1/|\si|), \quad |\nabla V|=O(1/|\si|^2), \quad |\si|\to\infty.
\]
Thus the first Green's identity implies the following relations for solutions of the homogeneous problem \rf{V0}-\rf{V0asy}, \rf{sys72}:
\begin{align*}
\int_{\Omh} |\nabla V|^2 \, \D \si =\int_{\de \Omh} \frac{\de V_-}{\de \n} \,V\, \D S ,
\qquad
\int_{\mathbb{R}^3 \setminus \Omh} |\nabla V|^2 \, \D \si  = \int_{\de \Omh} \frac{\de V_+}{\de \n} \,V\, \D S .
\end{align*}
We subtract the first equality multiplied by $\ep$ from the second one. Then homogeneous relation \rf{sys72} implies that
$V_\pm $ are constants. These constants are equal due to \rf{V0asy}, and therefore they are zeros due to \rf{Vasy}. Thus, uniqueness is established.
\end{proof}

Since the system of equations \rf{sys1} is represented in \rf{sysC} in the form of a power series in $a$ whose leading term is uniquely invertible due to Lemma \ref{lemma2}, solution of \rf{sys1} has the same property. That is, the following statement holds.
\begin{lemma}
 Solution $\wh\in H^{2}(\Omh),~ \walpha\in \Ho(\de\Omh)$ of the system of integral equations \rf{sys1} has the following power expansion for small values of $a$:
 \begin{align}
 \label{wh}
 \wh &= \wh_0 + \wh_1 a + \wh_2 a^2 + \ldots,\\[2mm]
 \walpha &= \walpha_0 + \walpha_1 a + \walpha_2 a^2 + \ldots
 \label{ah}
\end{align}
\label{lemma3}
\end{lemma}

An immediate corollary of the Lemma \ref{lemma3} is a similar statement about the solution of the system \rf{w}-\rf{alpha}:
\begin{lemma}\label{l10}
 The system of integral equations \rf{w}-\rf{alpha} has a unique solution $w\in H^{2}(\Om),~ \alpha\in \Ho (\de\Om)$ that can be represented in the form of a power series in $a$ whose coefficients depend smoothly on $a$:
 \begin{align}
 w &= w_0 + w_1 a + w_2 a^2 + \ldots,\quad w_n (\x,a) = \wh_n (a\x),\\[2mm]
 \alpha &= \alpha_0 + \alpha_1 a + \alpha_2 a^2 + \ldots,
 \quad \alpha_n (\x,a) = \walpha_n(a\x),
\end{align}
where series converge in the spaces $H^{2}(\Om),~ \Ho(\de\Om)$, respectively, when $a$ is sufficiently small.
\end{lemma}

\subsection*{Proof of Theorem \ref{t2}}

\begin{proof}
Without loss of generality, we can assume that $\Omh\subset B_{R/2}$ since this can be achieved by rescaling $a$.

We look for the solution of \rf{eq1}-\rf{bc1c} in the form \rf{ut}. This reduces operator $\Nm_{a,\varepsilon}-\Nm_{0,\varepsilon}$ to the map \rf{dn}. Formula \rf{sol} with $w, \alpha$ defined in Lemma \ref{l10} provides solution $\ut$ of \rf{eq2a}-\rf{bc1aa} for $a\ll  1$. We expand $G(\x,\y)$ in \rf{sol} in the Taylor series in $\y$. Due to \rf{G}, the latter expansion in $C^\infty$ when $\x$ is in a neighborhood of $\de B_R$ and $\y\in\Om$ has the form
 \[
 G(\x,\y)\sim \sum_{j=0}^\infty P_j(\x,\y),  \quad |\x|>R/2,~~|\y|\ll1,
 \]
 where $P_j$ are homogeneous polynomials in $\y$ of order $j$ with infinitely smooth in $\x$ coefficients. In particular,
 \begin{align}
 P_0 (\x,\y) &= \frac{\kp}{4\pi} \left[ y_0 (\kp |\x|) - j_0 (\kp |\x|) \frac{y_0 (\kp R)}{j_0 (\kp R)}\right], \\[2mm]
 P_1 (\x,\y) &= \frac{\kp^2}{4\pi} \left[ y_1 (\kp |\x|) - j_1 (\kp |\x|) \frac{y_1 (\kp R)}{j_1 (\kp R)}\right]
 \dotx \cdot \,\y, \quad \dotx \,= \x/|\x|.
 \end{align}
 Here $j_n(z), y_n(z)$ are the spherical Bessel functions of the first and second kind, respectively.
Hence, \rf{sol} allow us to represent solution $\ut$ of \rf{eq2a}-\rf{bc1aa} for $a\ll  1$ and $|\x|>R/2$ in the form
\[
\ut= \sum_{n=0}^N \sum_{j=0}^\infty  a^n \int_\Om  P_j(\x,\y) w_n(\y,a)d\y+\sum_{n=0}^N \sum_{j=0}^\infty  a^n \int_{\de\Om} P_j(\x,\y)\alpha_n(\y,a)dS_\y + \O \left(a^{N+1}\right),
\]
where the interior series converge in $C^\infty$. Since polynomials $P_j$ are homogeneous in $y$,
the substitution $\y\to a\et$ implies
\[
\ut\sim \sum_{j,n=0}^\infty a^{3+j+n}\int_{\Omh} P_j(\x,\et) \wh_n(\et)d\et+\sum_{j,n=0}^\infty a^{2+j+n}\int_{\de\Omh} P_j(\x,\et) \walpha_n(\et)dS_\et.
\]
Hence
\begin{align}
 \frac{\de\ut}{\de\n}\sim\sum_{n=2}^\infty a^n\varphi_n (\x),  \quad \x \in \de B_R,
 \label{dut}
\end{align}
where
\begin{align}
\label{int2}
\varphi_2 &=d\int_{\de\Omh}\walpha_0(\et)\D S_\et, \\[2mm]
\varphi_3 &=d\int_{\de\Omh}\walpha_1(\et)\D S_\et
+d_1\int_{\de\Omh}(\dotx \cdot \et)\,\walpha_0(\et)\D S_\et+d\int_{\Omh}\wh_0(\et)\D \et.
\label{int}
\end{align}
Here
\begin{align}
 d &= d(R) = \left. P_0^\prime \right|_{|\x|=R} = -\frac{\kp^2}{4\pi} \left[ y_1 (\kp R) - j_1 (\kp R) \frac{y_0 (\kp R)}{j_0 (\kp R)}\right] =  \frac{1}{4\pi R^2 j_0 (\kp R)}, \\[2mm]
 d_1 &= d_1 (R) = \left. P_1^\prime \right|_{|\x|=R} = -\frac{\kp^3}{4\pi} \bigl[ y_2 (\kp R) j_1 (\kp R) - j_2 (\kp R) y_1 (\kp R)\bigr] \frac{1}{j_1 (\kp R)} =  \frac{\kp}{4\pi R^2 j_1 (\kp R)},
 \label{dut1}
\end{align}
where we have used the property of the cross-product of the spherical Bessel functions \cite{Olver:10}
\begin{align}
 j_{n+1} (z) y_n (z) - j_n (z) y_{n+1}(z) = \frac{1}{z^2}.
 \label{}
\end{align}
The first statements of Theorem \ref{t2} will follow from \rf{ut}, \rf{dut}-\rf{int} if we show that $\phi_2=0$. We also need to specify the right-hand side in \rf{int} to prove the second statement of the theorem.

We evaluate integral in \rf{int2}.  From \rf{approx0} it follows that $\walpha_0$ satisfies \rf{ttt} with $f = h_0 = \varkappa \n \cdot {\bm c}$. Hence Lemma \ref{l3} with $\alpha = \alpha_0$ yields
\begin{align}\nonu
 \walpha_0 &= \frac{\de V_+}{\de \n} - \frac{\de V_-}{\de \n}
 = \frac{\de V_+}{\de \n} - \ep\, \frac{\de V_-}{\de \n} + \ep\, \frac{\de V_-}{\de \n} - \frac{\de V_-}{\de \n} \\[2mm]
 &= (1+\ep)h_0 + (\ep -1)\, \frac{\de V_-}{\de \n}
 = (\ep -1) (\n \cdot {\bm c}) + (\ep -1)\, \frac{\de V_-}{\de \n}, \quad \si\in \de \Omh.
 \label{alpha0}
\end{align}
Thus,
\begin{align}
 \phi_2 = d\int_{\de \Omh} \walpha_0\, \D S = (\ep-1) d \int_{\de \Omh} (\n \cdot {\bm c})\, \D S + (\ep - 1)d \int_{\de \Omh} \frac{\de V_-}{\de \n} \, \D S.
 \label{fi1}
\end{align}
The latter integral vanishes due to Green's first identity,  and the divergence theorem implies that
\begin{align}
 \phi_2 = (\ep-1) d\int_{\de \Omh} (\n \cdot {\bm c})\, \D S
 = (\ep -1) d\int_{\Omh} \left( \nabla \cdot {\bm c} \right) \D V = 0.
 \label{fi2}
\end{align}
Hence, the first statement of the theorem is proven.

Next, we study \rf{int}. Since $\phi_3=\P_3(\varepsilon)\psi$ where $\P_j$ are defined in \rf{et2}, we need formula \rf{int} only for $\psi = \hpsi_i = \E^{-\I (\bk - \mb_i)\cdot \x}$, and therefore we will use $u(\z) = 1$,
${\bm c} = -\I (\bk - \mb_i)$, $\Q$ is the negative dyadic square of the vector $(\bk - \mb_i)$, and $\tr \,\Q = -|\bk - \mb_i|^2 = -|\bk|^2$.  The same argument used to evaluate the integral of $\walpha_0$ can be applied to evaluate
the integral of $\walpha_1$. From \rf{approx1} it follows that $\walpha_1$ satisfies \rf{ttt} with $f = h_1 + \varkappa \what{\G}^\prime_1 \wh_0 = \varkappa {\Q \si} \cdot {\n}+ \varkappa \what{\G}^\prime_1 \wh_0 $.  Hence Lemma \ref{l3} with $\alpha = \alpha_1$ yields
\begin{align}
 \walpha_1 &= \frac{\de V_+}{\de \n} - \frac{\de V_-}{\de \n}
 = \frac{\de V_+}{\de \n} - \ep\, \frac{\de V_-}{\de \n} + \ep\, \frac{\de V_-}{\de \n} - \frac{\de V_-}{\de \n} \nonu \\[2mm]
 &= (1+\ep)\left( h_1 + \varkappa \what{\G}^\prime_1 \wh_0 \right) + (\ep -1)\, \frac{\de V_-}{\de \n}, \quad \si\in \de \Omh.
 \label{alpha1}
\end{align}
Since the surface integral over $\de \Omh$ of the latter term in \rf{alpha1} vanishes, from formula \rf{u1exp} for $h_1$ and expression for $\varkappa$ from Lemma \ref{lemma1} it follows that
\begin{align*}
\int_{\de\Omh}\walpha_1(\et)\,\D S_\et &=
 (\ep-1)\int_{\de \Omh} (\Q \si \cdot \n) \, \D S + (\ep-1) \int_{\de \Omh} \what{\G}^\prime_1 \wh_0 \, \D S. 
\end{align*}
Using formulas \rf{Gpp}, \rf{T0h}  for $\what{\G}^\prime_1$, formula \rf{w0} for $\wh_0$ and the divergence theorem we obtain
\begin{align}
\int_{\de\Omh}\walpha_1(\et)\,\D S_\et &=
  (\ep-1) \int_{\Omh} \nabla \cdot \left(\Q \si\right) \, \D V + (\ep-1)  (\kp^2 - \km^2) \int_{\de \Omh} \int_{\Omh} \what{T}_0 (\si, \et) \, \D V\, \D S.
 \label{alpha1int}
\end{align}
We change the order of integration in the last integral. Since $\int_{\de \Omh} \what{T}_0 (\si, \et) \,\D S = 1, \; \et \in \Omh$, the last term in the right-hand side of \rf{alpha1int} equals $(\ep -1)  (\kp^2 - \km^2)  |\Omh |$. The first term equals $(1-\ep) |\bk|^2|\Omh|$. Hence,
\begin{align}
\int_{\de\Omh}\walpha_1(\et)\,\D S_\et &= (1-\ep) \left( |\bk|^2 -   ( \kp^2 -\km^2) \right)|\Omh|.
 \label{int_a1}
\end{align}
Due to \rf{w0} we have
\begin{align}
 \int_{\Omh}\wh_0(\et)\,\D \et = (\kp^2 - \km^2) \int_{\Omh} \D \et = (\kp^2 - \km^2) |\Omh |.
 \label{int_w0}
\end{align}

It remains to evaluate the middle term in \rf{int}:
\begin{align*}
 \J := d_1\int_{\de\Omh}(\dotx \cdot \et)\,\walpha_0(\et)\D S_\et.
\end{align*}
This integral is related to the expansion of the simple layer potential \rf{vh} in the form of  a power series in $|\si|^{-1}$ when $|\si| \gg |\et|$:
 \begin{align}
  V&= -\int_{\de \Omh}\frac{1}{4\pi |\si - \et|}\, \walpha_0 (\et)\, \D S_{\et}
  = -\frac{1}{4\pi | \si|}\int_{\de \Omh} \walpha_0 (\et)\left( 1 + \frac{\si \cdot \et}{|\si|^2} + \ldots \right)\, \D S_{\et} \nonu \\[2mm]
  &= -\frac{1}{4\pi} \left( \frac{p_0}{|\si|} + \frac{\si \cdot \m}{|\si|^3} + \ldots \right),
  \label{Vexp}
 \end{align}
where
\begin{align}
 p_0 &= \int_{\de \Omh} \walpha_0 (\et)\, \D S_{\et}, \quad
 \m = \int_{\de \Omh} \et\,\walpha_0 (\et)\, \D S_{\et}.
 \label{p0p}
\end{align}
Due to \rf{fi1}-\rf{fi2}, $p_0 =0$. Hence \rf{Vexp},  \rf{p0p} provide \rf{ppp}  and
\begin{align}
 \J = d_1\int_{\de\Omh}(\dotx \cdot \et)\,\walpha_0(\et)\D S_\et =
  d_1 {\bm p}\cdot \dotx,
 \label{J}
\end{align}
where $\m=\m_i$ corresponds to the external field $\E^{\I(\bk - {\bm m}_i) \cdot \x}$.
Combining \rf{int}, \rf{int_a1}, \rf{int_w0} and \rf{J} we obtain:
\begin{align}
 \phi_3 &=\P_3(\varepsilon)\hpsi_i= (1-\ep) d \left(|\bk|^2  -  (\kp^2 - \km^2)  \right) |\Omh| + d (\kp^2 - \km^2) |\Omh |
 + d_1  {\bm p}_i\cdot \dotx \nonu \\[2mm]
 &=(1-\ep) d  |\bk|^2 |\Omh|  + \ep d (\kp^2 - \km^2) |\Omh | + d_1  {\bm p}_i\cdot \dotx  = d |\bk|^2 |\Omh|  - \ep d  \km^2 |\Omh | + d_1  {\bm p}_i\cdot \dotx+O(\varepsilon) \nonu\\[2mm]
 &=  d   |\bk|^2  \left(1 -\frac{\gi}{\ge} \right) |\Omh | + d_1  {\bm p}_i\cdot \dotx +O(\varepsilon),  \quad \ei=(\kp -|\bk|)/|\bk|.
 \label{fp}
\end{align}
We used here that $\kp = |\bk| + \O(\varepsilon)$ and $\ep \km^2 = \frac{\gi}{\ge}\, |\bk|^2 + \O(\varepsilon)$. Formula \rf{fp} together with \rf{et2} and \rf{chi1} implies
\[
(\Nm_{a,\varepsilon} -\Nm_{0,\varepsilon}) \psi_{i}=d   |\bk|^2  \left(1 -\frac{\gi}{\ge} \right) |\Om |
- \I d_1 |\bk| |\Om|  {\bm \chi_i}\cdot \dotx +O(\varepsilon a^3+a^4)
\]
Hence, the elements $M_{ij},~ 0\leq i,j\leq n,$ of matrix $\M$ are
\begin{align}
M_{i,j}&=\left( (\Nm_{a,\varepsilon} -\Nm_{0,\varepsilon}) \psi_{i}, \psi_{j} \right)=d   |\bk|^2  \left(1 -\frac{\gi}{\ge} \right) |\Om |\int_{\de B_R} \E^{\I (\bk-\mb_j) \cdot \x}\, \D S \nonu \\[2mm]
 &- \I d_1 |\bk| |\Om|  {\bm \chi_i}\cdot \int_{\de B_R}\dotx\, \E^{\I (\bk-\mb_j) \cdot \x}\, \D S+O(\varepsilon a^3+a^4).
\label{mij}
\end{align}

Let us evaluate the integrals above. We have
\[
\int_{\de B_R} \E^{\I (\bk-\mb_j) \cdot \x}\, \D S=4\pi R^2 j_0 (|\bk-\mb_j| R)=4\pi R^2 j_0 (|\bk| R)=4\pi R^2 j_0 (\kp R)+O(\varepsilon)=\frac{1}{d}+O(\varepsilon).
\]
\begin{align}
& \int_{\de B_R} \dotx\,   \E^{\I (\bk-\mb_j) \cdot \x}\, \D S =R^{-1}\int_{\de B_R} \x  \E^{\I (\bk-\mb_j) \cdot \x}\, \D S =-\I R^{-1}\nabla_\bk \int_{\de B_R}  \E^{\I (\bk-\mb_j) \cdot \x}\, \D S \nonu \\[2mm]
&=-4\pi\I R\nabla_\bk j_0 (|\bk-\mb_j| R)=4\pi\I R^2 j_1 (|\bk-\mb_j| R)\frac{\bk-\mb_j}{|\bk-\mb_j|}=\I|\bk|\frac{\bk-\mb_j} {d_1|\bk-\mb_j|}+\O \left(\varepsilon \right).
\end{align}
Substitution of the values of the integrals into \rf{mij} leads to \rf{Mat2}, and this completes the proof of the theorem.
\end{proof}

\section*{Acknowledgment}

The work of B. Vainberg was supported by the Simons Foundation grant 527180.


\end{document}